\documentclass[a4paper, 10pt]{article}

\usepackage{amsthm, amsmath, amssymb, amsfonts}
\usepackage[utf8]{inputenc}
\usepackage[english]{babel}
\usepackage{url}
\usepackage{mathtools}
\usepackage{romannum}

\numberwithin{figure}{section}
\numberwithin{table}{section}
\numberwithin{equation}{section}

\newenvironment{abstr}[1]{ \vspace{.05in}\footnotesize
       \parindent .2in
         {\upshape\bfseries #1. }\ignorespaces}{\par\vspace{.1in}}
\newenvironment{Abstract}{\begin{abstr}{Abstract}}{\end{abstr}}
\newenvironment{keywords}{\begin{abstr}{Key words}}{\end{abstr}}
\newenvironment{AMS}{\begin{abstr}{AMS subject classifications}}{\end{abstr}}

\newtheorem{theorem}{Theorem}[section]
\newtheorem{lemma}[theorem]{Lemma}

\newtheorem{proposition}[theorem]{Proposition}
\newtheorem{definition}[theorem]{Definition}
\newtheorem{assumption}[theorem]{Assumption}
\newtheorem{remark}[theorem]{Remark}
\newtheorem{notation}[theorem]{Notation}
\newtheorem{conclusion}[theorem]{Conclusion}

\newenvironment{Proof}[1][Proof]{\begin{proof}[#1]}{\end{proof}}

\DeclareMathOperator*{\TO}{\longrightarrow}
\DeclareMathOperator{\supp}{supp}
\DeclareMathOperator{\diam}{diam}
\DeclareMathOperator{\Div}{div}
\DeclareMathOperator{\curl}{curl}

\renewcommand{\Re}{\operatorname{Re}}
\renewcommand{\Im}{\operatorname{Im}}
\DeclareMathOperator{\Res}{Res}
\newcommand{\halb}{\frac 12}
\DeclareMathOperator*{\wto}{\rightharpoonup}
\newcommand{\twosc}{\stackrel{2}{\wto}}

\newcommand{\nz}{\mathbb{N}}       
\newcommand{\rz}{\mathbb{R}}       
\newcommand{\cz}{\mathbb{C}}       
\newcommand{\pz}{\mathbb{P}}
\newcommand{\de}{\delta}

\newcommand{\ep}{\varepsilon}

\newcommand{\om}{\omega}
\newcommand{\Om}{\Omega}

\newcommand{\si}{\sigma}

\newcommand\Va{\mathbf{a}}
\newcommand\Vb{\mathbf{b}}

\newcommand\Ve{\mathbf{e}}
\newcommand\Vf{\mathbf{f}}

\newcommand\Vj{\mathbf{j}}

\newcommand\Vn{\mathbf{n}}

\newcommand\Vv{\mathbf{v}}
\newcommand\Vu{\mathbf{u}}
\newcommand\Vw{\mathbf{w}}

\newcommand\Vz{\mathbf{z}}

\newcommand\VE{\mathbf{E}}
\newcommand\VF{\mathbf{F}}
\newcommand\VH{\mathbf{H}}

\newcommand\VK{\mathbf{K}}

\newcommand\VR{\mathbf{R}}

\newcommand\VV{\mathbf{V}}
\newcommand\VW{\mathbf{W}}
\newcommand\Vpsi{\boldsymbol{\psi}}

\newcommand\CA{\mathcal{A}}
\newcommand\CB{\mathcal{B}}

\newcommand\CE{\mathcal{E}}

\newcommand\CT{\mathcal{T}}

\newcommand{\Hbf}{\mathbf{H}}

\usepackage[textsize=footnotesize]{todonotes}
\usepackage[xcolor=dvipdf]{changes}
\setremarkmarkup{\todo[color=Changes@Color#1!40]{\textcolor{black}{#1:{#2}}}}

\allowdisplaybreaks[4]

\title{A new Heterogeneous Multiscale Method for time-harmonic Maxwell's equations based on divergence-regularization\thanks{This work was supported by the Deutsche Forschungsgemeinschaft (DFG) in the project "Wellenausbreitung in periodischen Strukturen und Mechanismen negativer Brechung"}}
\author{Patrick Henning\footnotemark[2] \and Mario Ohlberger\footnotemark[2] \and Barbara Verf\"urth\footnotemark[2]}
\date{}

\begin{document}
\pagenumbering{arabic}

\maketitle
\renewcommand{\thefootnote}{\fnsymbol{footnote}}
\footnotetext[2]{Institut f\"ur Numerische und Angewandte Mathematik, Westf\"alische Wilhelms-Uni\-ver\-si\-t\"at M\"unster, Einsteinstr. 62, D-48149 M\"unster}
\renewcommand{\thefootnote}{\arabic{footnote}}

\begin{Abstract}
In this paper, we suggest a new heterogeneous multiscale method (HMM) for the time-harmonic Maxwell equations in locally periodic media. 
The method is constructed by using a divergence-regularization in one of the cell problems. This allows us to introduce fine-scale correctors 
that are not subject to a cumbersome divergence-free constraint and which can hence easily be implemented.
To analyze the method, we first revisit classical homogenization theory for time-harmonic Maxwell equations and derive a new 
homogenization result that makes use of the divergence-regularization in the two-scale homogenized equation.
We then show that the HMM is equivalent to a discretization of this equation. In particular, writing both problems in a fully coupled two-scale formulation is the crucial starting point for a corresponding numerical analysis of the method. 
With this approach we are able to prove rigorous 
a priori error estimates in the $\Hbf(\curl)$- and the $H^{-1}$-norm 
and we derive 
reliable and efficient localized residual-based a posteriori error estimates.
\end{Abstract}



\section{Introduction}
The behavior and propagation of electromagnetic fields is studied in many physical applications, for instance in the large area of wave optics. 
Periodic and locally periodic materials are considered with growing interest, for example in the application of photonic crystals (see \cite{phC} for an introduction), as they can show unusual behavior, such as photonic band gaps and even negative refraction (see e.g.\ \cite{EP04, LJJP02, PE03lhm, SPW04meta}).
However, a thorough mathematical understanding of these phenomena is still lacking. Therefore, one major goal is to develop efficient numerical schemes  to simulate wave propagation in periodic materials and to rigorously analyze the new algorithms and the errors they introduce. 

Maxwell's equations in a linear conductive medium, subject to Ohm's law, are given by
\begin{align*}
\curl \hat{\VE}(x, t) +\mu(x) \partial_t \hat{\VH}(x,t) &=0,  &\Div(\mu(x) \hat{\VH}(x,t)) &=0,\\
\curl \hat{\VH}(x,t) -\ep(x) \partial_t \hat{\VE}(x,t) -\si(x) \hat{\VE}(x,t) &=\hat{\Vj}(x,t), & \Div(\ep(x)\hat{\VE}(x,t))&=\hat{\rho}(x,t).
\end{align*}
Here, the electric field $\hat{\VE}$ and the magnetic field $\hat{\VH}$ are the unknowns, the current density $\hat{\Vj}$ and the charge density $\hat{\rho}$ are the sources, and the permittivity $\ep$, the permeability $\mu$, and the conductivity $\si$ are material parameters. The time-harmonic Maxwell equations can be obtained by assuming periodicity in time, i.e.\ for $\psi=\VE,\VH,\Vj,\rho$ we can make
the ansatz $\hat{\psi}(x,t) =\Re(\psi(x)e^{i\om t})$ with a complex-valued function $\psi$ and a temporal frequency $\om\neq 0$. Inserting this ansatz in the original equations yields the time-harmonic system
\begin{align}
\label{eq:curlE}
\curl(\mu^{-1} \curl \VE)+(i\om \si-\om^2\ep)\VE &=-i\om \Vj,\\
\VH&=i\om^{-1}\mu^{-1}\curl\VE.
\end{align}

In this paper, we will focus on the curl-curl-problem \eqref{eq:curlE} on a bounded domain and with locally periodic coefficients. More precisely, let $\Om\subset \rz^3$ be a bounded domain with outer normal $\Vn$ on $\partial\Om$ and we seek $\VE_\de: \Om\to \cz^3$ with
\begin{align}
\label{eq:Edestrong}
\curl(\mu_\de^{-1}\curl \VE_\de)-\kappa_\de\VE_\de &=\Vf &&\text{in } \Om,\\
\label{eq:Edebdry}
\VE_\de\times \Vn &=0 &&\text{on }\partial\Om.
\end{align}
We will assume that the real-valued parameter $\mu^{-1}_\de$ and the complex-valued parameter $\kappa_\de$ are locally periodic with periodicity length $\de$, where $\de$ is very small compared to $\Om$.
The boundary condition \eqref{eq:Edebdry} models the case where $\Om$ is surrounded by a so called perfect electric conductor. We refer to \cite{Monk,Schoeberl}, or \cite{Zagl} for a detailed motivation and further applications.

Since a numerical treatment of \eqref{eq:Edestrong}--\eqref{eq:Edebdry} requires discretizations with mesh sizes $h<\delta\ll 1$, corresponding computations can easily exceed today's available computer resources if tackled with a standard approach. In order to make the problem numerically solvable, so called multiscale methods can be applied. One class of multiscale methods that has been proved to be very efficient for scale-separated problems with local periodicity (or mild heterogeneities) is the family of Heterogeneous Multiscale Methods (HMM) introduced by E and Engquist \cite{Eng,Eng1}. HMM approaches exploit structural invariants in the coefficients to solve local sample problems that allow to extract representative features and to approximate the $\delta$-dependent multiscale solution with a computational complexity that is however independent of $\delta$. With this strategy the problem becomes solvable even for arbitrarily small values of $\delta$. First analytical results concerning the approximation properties of the HMM for elliptic problems have been derived in \cite{Abdulle, EMZ05, Gloria:2006, Ohl}. In this contribution we formulate and analyze a new HMM for solving the curl-curl-problem \eqref{eq:Edestrong}--\eqref{eq:Edebdry}.

Concerning wave propagation in general, the HMM and related multiscale methods for wave equations have been studied by Abdulle et.al.\ \cite{Grote, Grote1} and Engquist,
Runborg 
et.al.\ \cite{AR14wave,  EHR11wave, EHR12wavelong, ER03wavehighfreq}. An HMM for the Helmholtz equation has been suggested in \cite{CS14HMMhelmh}. Furthermore, some methods based on asymptotic expansions have been suggested for Maxwell's equations, see e.g.\ \cite{CZAL, ZCW}. 

The new contribution of this article is the first formulation of a Heterogeneous Multiscale Method for the time-harmonic Maxwell equations and its comprehensive numerical analysis in terms of a priori and a posteriori error estimates.
The error analysis can serve as a starting point for a locally adaptive version of the described HMM and the scheme itself might be applicable to other related problems after slight modifications.

The idea of the HMM is to adapt the (analytical) homogenization procedure to the numerical scheme. Therefore, we will first have a look at the homogenization of the time-harmonic Maxwell equations. Combining results by Wellander et.al.\ \cite{SvanWell,Well1,Well2} and Visintin \cite{Visintin}, we derive a new two-scale equation for time-harmonic Maxwell's equations. One essential step in the homogenization procedure is a divergence-regularization in order to incorporate a divergence-free constraint imposed on the corrector of the curl into the equation. This regularization also is an essential ingredient in the formulation of the new HMM. 
We will then adopt the view of the HMM as a direct discretization of the derived two-scale equation. This reformulation builds the crucial ingredient for an a posteriori analysis. It has been first developed in \cite{Ohl} and has then been adopted to other problems, as perforated domains \cite{HO09perf} or advection-diffusion problems \cite{HO10advec}, for instance.
There have been several contributions on the numerical analysis for time-harmonic Maxwell's equations (see the excellent book \cite{Monk} for an overview, \cite{Zagl} for higher order finite elements, and \cite{Sch1, BHHW} for a posteriori analysis). In the analysis of our HMM we need to combine these two approaches in a new way. Thereby we are able to obtain optimal error estimates.

The article is organized as follows: In Section \ref{sec:problem} we formulate the multiscale curl-curl-problem and give some properties of the solution. The problem is homogenized with the tool of two-scale convergence in Section \ref{sec:hom}. The homogenized formulation is the motivation and starting point for the formulation of the HMM in Section \ref{sec:hmm}. Error estimates for this method are given in Section \ref{sec:error}. All essential proofs are detailed in Section \ref{sec:errorproof}.

\section{Problem setting}
\label{sec:problem}
For the remainder of this article, let $\Om\subset \rz^3$ be a bounded, simply connected domain with connected piecewise polygonal Lipschitz boundary $\partial \Om$ and outer unit normal $\Vn$.
Throughout this paper, we use standard notation: By $W^{l,p}(\Om)$ we denote  the space of functions on $\Om$ with weak derivatives up to order $l$ belonging to $L^p(\Om)$ and we write $H^l(\Om):=W^{l,2}(\Om)$ for the scalar and $\mathbf{H}^l(\Om):=[H^l(\Om)]^3$ for the vector-valued case. Vector-valued functions are indicated by boldface letters and unless otherwise stated, all functions are complex-valued. The dot will denote a normal (real) scalar product, for a complex scalar product we will explicitly conjugate the second component by using $u^*$ as the conjugate complex of $u$.
For any domain $\omega\subset \mathbb{R}^3$, we introduce the spaces
\begin{equation*}
\begin{split}
\Hbf(\curl,\omega)&:=\{\Vu\in L^2(\omega; \cz^3)| \hspace{2pt} \curl \Vu\in L^2(\omega; \cz^3)\} \qquad \mbox{and}\\
\Hbf(\Div,\omega)&:=\{\Vu\in L^2(\omega; \cz^3)| \hspace{2pt} \Div \Vu\in L^2(\omega; \cz)\}.
\end{split}
\end{equation*}
For $\omega=\Om$ we write $\Hbf(\curl):=\Hbf(\curl,\Omega)$ and $\Hbf(\Div):=\Hbf(\Div,\Omega)$. These spaces are complex Hilbert spaces if endowed with the scalar products
\begin{equation*}
\begin{split}
(\Vu, \Vv)_{\Hbf(\curl)}&:=\int_\Om \curl \Vu\cdot \curl \Vv^*+\Vu\cdot \Vv^*\, dx,\\
(\Vu, \Vv)_{\Hbf(\Div)}&:=\int_\Om \Div \Vu\, \Div\Vv^*+\Vu\cdot\Vv^*\, dx.
\end{split}
\end{equation*}
With the help of a trace theorem (see \cite{Monk}), zero boundary values for functions in $\Hbf(\curl)$ can be defined as
$$\Hbf_0(\curl):=\{\Vv\in \Hbf(\curl)|\Vv\times \Vn=0\}.$$
For higher regularity, we define for $s\in \nz_0$ the space
$$\Hbf^s(\curl):=\{\Vu\in \Hbf(\curl)\, |\; \Vu\in \Hbf^s(\Om), \curl\Vu\in \Hbf^s(\Om)\}.$$
Observe that $\Hbf^0(\curl)=\Hbf(\curl)$.
Let $\Ve_i$ denote the $i$'th unit vector in $\rz^3$ (i.e. $(\Ve_i)_j = \delta_{ij}$ for $1\le i,j \le 3$). For the rest of the paper we write $Y:=[-\halb, \halb)^3$ to denote the 3-dimensional unit cube and we say that a function $v\in L^2_{\mbox{\tiny loc}}(\mathbb{R}^3)$ is $Y$-periodic if it fulfills $v(y)=v(y+\Ve_i)$ for all $i=1,2,3$ and almost every $y\in \mathbb{R}^3$. With that we denote $L_{\sharp}^2(Y):=\{ v \in L^2_{\mbox{\tiny loc}}(\mathbb{R}^3)| \hspace{2pt} v \mbox{ is $Y$-periodic}\}$. Analogously we indicate periodic function spaces by the subscript $\sharp$. For example, $H^1_\sharp(Y)$ is the space of periodic $H^1_{\mbox{\tiny loc}}(\mathbb{R}^3)$-functions and we furthermore define for $s\in\mathbb{N}$
$$H^s_{\sharp,0}(Y):=\left\{ \left.\phi \in H^s_\sharp(Y)\right|\int_Y \phi(y) \, dy =0\right\}.$$
By $L^p(\Om; X)$ we denote Bochner-Lebesgue spaces over the Banach space $X$ and we use the short notation $f(x,y):=f(x)(y)$ for $f\in L^p(\Om; X)$.

Using the above notation we make the following assumptions on the coefficients.
\begin{assumption}
\label{assumption-coefficients}
The (scalar) coefficient $\mu^{-1} \in C^0(\Om; L^\infty_\sharp(Y))$ is real-valued and $\kappa \in C^0(\Om; L^\infty_\sharp(Y;\mathbb{C}))$ is complex-valued. 
Let $\psi$ denote $\Re(\kappa)$, $-\Im(\kappa)$ or $\mu^{-1}$.
Then there exist $c_0,c_1\in\mathbb{R}$ such that
$0<c_0\leq \psi(x,y) \leq c_1$ for a.e. $x$ and $y$, such that $\psi(\cdot,\frac{\cdot}{\delta})$ is measurable for all $\delta>0$ and
\begin{equation}
\label{eq:assptstrong2sc}
\lim_{\de\to 0}\int_\Om \psi\Bigl(x, \frac{x}{\de}\Bigr)^2\, dx =\int_\Om\int_Y \psi(x,y)^2\, dydx.
\end{equation}
\end{assumption}

\begin{definition}[Weak solution]
\label{def:weakEde} 
Define $\kappa_\de(x):=\kappa(x, \frac{x}{\de})$, $\mu^{-1}_\de(x):=\mu^{-1}(x, \frac{x}{\de})$ and let Assumption \ref{assumption-coefficients} be fulfilled. 
 Let $\Vf\in \Hbf(\Div)$ with $\Div \Vf=0$. We call 
$\VE_\de\in \Hbf_0(\curl)$ a weak solution if, for all $\Vpsi\in \Hbf_0(\curl)$, it fulfills
\begin{equation}
\label{eq:weakEde}
\begin{split}
&\int_\Om \mu^{-1}_\de(x)\curl \VE_\de(x) \cdot \curl \Vpsi^*(x)-\kappa_\de\VE_\de(x)\cdot \Vpsi^*(x)\, dx
=\int_\Om \Vf(x)\cdot \Vpsi^*(x)\, dx.
\end{split}
\end{equation}
\end{definition} 

For fixed $\de$, there is a unique solution to \eqref{eq:weakEde}, which can be seen using the Lax-Milgram-Babu{\v{s}}ka theorem, \cite{Bab70}: Clearly, the right-hand side is a member of the dual space and the left-hand side gives a continuous sesquilinear form. Since $\Im\kappa$ is bounded away from zero, we also get the coercivity estimate $|B_\de(\Vu, \Vu)|\geq C\|\Vu\|^2_{\Hbf(\curl)}$ with a $\de$-independent constant. See \cite{Zagl} for the case of constant coefficients and \cite{FeRa} for the general computation. Hence, we also have the uniform estimate $\|\VE_\de\|_{\Hbf(\curl)}\leq C\|\Vf\|_{L^2}$ with $C=C(c_0,c_1,\Omega)$.

In general, solutions to curl-curl-problems do not admit more than $H^{1/2}$-regularity and may have singularities near re-entrant corners of the domain, see \cite{CD}. However, if $\Om$ is convex and if $\mu^{-1}, \kappa\in W^{1,\infty}(\Om\times Y)$, i.e.\ the coefficients are globally Lipschitz, the weak solution to \eqref{eq:weakEde} fulfills $\VE_\de\in \Hbf^1(\curl)$, see \cite{Schoeberl}.

\section{Homogenization}
\label{sec:hom}

As the periodicity length $\de$ is assumed to be very small in comparison to $\Om$, one can reduce the complexity of the problem by considering the limit $\de \to 0$. This process is called homogenization and can be performed with the tool of two-scale convergence \cite{Alla}. Since the two-scale equation and the homogenized equation are the starting point for the construction and analysis of the numerical multiscale method later on, we present the essential steps in this section.

\subsection{Two-scale convergence}

Two-scale convergence is defined as (cf. \cite{Alla,LNW}):
\begin{definition}[Two-scale convergence] 
\label{def:zweiskKonv}
A sequence $(u_\de)_{\de >0}\subset L^2(\Om)$ two-scale converges to a function $u_0\in L^2(\Om\times Y)$ (short form: $u_\de\twosc u_0$) if it fulfills
$$\lim_{\de \to 0}\int_\Om u_\de(x)\psi\left(x,\frac{x}{\de}\right)\, dx =\int_\Om\int_Y u_0(x,y)\psi(x,y)\, dydx\qquad \forall \psi \in L^2(\Om; C^0_\sharp(Y)).$$
\end{definition}

For more information on two-scale convergence, for example the definition of strong two-scale convergence, and compactness in $L^2$ and $H^1$, we refer the reader to  \cite{Alla,LNW}, and the lecture script \cite{Henn}.

For time-harmonic Maxwell's equations we need a two-scale convergence result for bounded sequences in $\Hbf(\curl)$. As $\Hbf(\curl)$ is not compactly embedded in $L^2$ (in contrast to $H^1$), the two-scale limit in $L^2$ will not coincide with the weak limit thus making additional considerations necessary (cf.\ \cite[p.\ 135]{Visintin}). We present two possible approaches, which we will combine in our analysis later on.
With the help of the $L^2$ compactness theorem and integration by parts, the following characterization of two-scale limits can be derived (see \cite{SvanWell,Well1,Well2}):

\begin{proposition}
\label{thm:Well}
Let $(\Vu_\de)_{\de>0}\subset \Hbf(\curl)$ be a bounded sequence. Then there exists a subsequence and functions $\Vu_0\in L^2(\Om\times Y)$, $\tilde{\Vu}_1\in L^2(\Om; \Hbf_{\sharp}(\curl,Y))$ and $\phi\in L^2(\Om; H^1_{\sharp, 0}(Y))$ such that
\begin{enumerate}
\item $\Vu_\de\twosc \Vu_0$ with $\Vu_0(x)=\Vu(x)+\nabla_y \phi(x,y)$ and $\Vu=\int_Y\Vu_0(\cdot,y)\, dy\in \Hbf(\curl)$,
\item $\curl\Vu_\de\twosc \curl_x \Vu_0+\curl_y\tilde{\Vu}_1$ and $\curl\Vu_\de\wto\curl\Vu$ in $L^2(\Om;\mathbb{C}^3)$.
\end{enumerate}
\end{proposition}

Using a technique called periodic unfolding, one can obtain a characterization for the curl which resembles the one in the $H^1$ compactness theorem (see \cite{Visintin}):

\begin{proposition}
\label{thm:Visintin}
Let $(\Vu_\de)_{\de>0}\subset \VH(\curl)$ be a bounded sequence such that $\Vu_\de\twosc \Vu_0$ in $L^2(\Om\times Y)$. Then there exists a subsequence and $\Vu_1\in L^2(\Om; \Hbf^1_{\sharp,0}(Y) )$ with $\Div_y \Vu_1=0$ a.e.\ in $\Om\times Y$ such that 
$$\curl \Vu_\de \twosc \curl \Vu+\curl_y \Vu_1 \text{ with } \Vu:=\int_Y\Vu_0(\cdot,y)\, dy.$$
Furthermore, it holds $\Vu_0\in L^2(\Om; \Hbf_\sharp(\curl, Y))$ and we have for the whole sequence $\de\curl \Vu_\de\twosc\curl_y \Vu_0=0$.
\end{proposition}
We will combine this theorem with the first point of Theorem \ref{thm:Well}. Note that the condition $\Div_y \Vu_1=0$ can be seen as  a kind of gauging condition. It will be important for the homogenization of our curl-curl-problem, namely this condition will lead to the uniqueness of the two-scale solution. 

\subsection{Homogenization for time-harmonic Maxwell's equations}

In this section we present new
homogenization results for the time-harmonic Maxwell equations 
in a two-scale formulation,
a formulation with cell problems and macroscopic equations, and a corrector result. We emphasize that although Maxwell's equations and curl-curl-problems have been homogenized in the literature (\cite{Amirat,Well1,Well2} only to name a few), the focus has always been on macroscopic (homogenized) problems as \eqref{eq:macroE}, but 
not
on two-scale limit equations. 

\begin{theorem}[Two-scale equation]
\label{thm:twosc}
Under the same assumptions as in Definition \ref{def:weakEde}, let $\VE_\de\in \Hbf_0(\curl)$ be a solution of \eqref{eq:weakEde}. Then there exists a solution triple 
$(\VE, \VK_1, K_2)$ of functions $\VE\in \Hbf_0(\curl)$, $\VK_1\in L^2(\Om; \Hbf^1_{\sharp,0}(Y))$ with $\Div_y\VK_1=0$ a.e., and $K_2\in L^2(\Om; H^1_{\sharp,0}(Y))$ such that
$$\VE_\de\wto \VE \text{ in }\Hbf_0(\curl),\quad \VE_\de\twosc\VE+\nabla_yK_2, \quad\curl\VE_\de\twosc\curl\VE+\curl_y\VK_1.$$ 
Considered in $\Hbf_0(\curl)\times L^2(\Om; \Hbf^1_{\sharp,0}(Y))\times  L^2(\Om; H^1_{\sharp,0}(Y))$, the triple $(\VE, \VK_1, K_2)$ is the unique solution of
\begin{equation}
\label{eq:twoscE}
\begin{split}
&\!\!\int_\Om\int_Y\mu^{-1}(x, y)(\curl\VE(x)+\curl_y\VK_1(x,y))\cdot(\curl\Vpsi^*(x)+\curl_y\Vpsi_1^*(x,y))\\
&\qquad\qquad +\Div_y \VK_1(x,y)\Div_y \Vpsi_1^*(x,y)\\
&\qquad \qquad-\kappa(x, y)(\VE(x)+\nabla_yK_2(x,y))\cdot(\Vpsi^*(x)+\nabla_y\psi_2^*(x,y))\, dy dx\\
&=\int_\Om \Vf(x)\cdot\Vpsi^*(x)\, dx \, 
\\&\qquad\forall\Vpsi\in \Hbf_0(\curl), \Vpsi_1\in L^2(\Om; \Hbf^1_{\sharp,0}(Y)), \psi_2\in L^2(\Om; H^1_{\sharp,0}(Y)).
\end{split}
\end{equation}
\end{theorem}

The proof is postponed to Section \ref{sec:errorproof}, but let us name the three important steps it consists of. First, using the two-scale convergence results mentioned above and inserting a special test function, we obtain a two-scale equation similar to \eqref{eq:twoscE}, but with divergence-free constraint. Second, we incorporate the divergence-free constraint into the equation (see also the remarks below). Third, the uniqueness of the two-scale solution is shown, which gives the convergence for the whole sequence of solutions $\VE_\de$.

\begin{remark}[Divergence-regularization]
\label{rem:divreg}
In order to determine $\VK_1$ in the two-scale equation, one 
has to find 
$\Vu \in \Hbf(\curl, Y)\cap \Hbf(\Div, Y)$ with $\Div \Vu =0$ a.e. in $Y$ and such that for all $\Vpsi\in \Hbf(\curl, Y)$
\begin{align}
\label{eq:curlcurl}
\int_Y \mu^{-1}\curl \Vu\cdot \curl\Vpsi^*\, dy &=0,
\end{align}
and with appropriate boundary conditions, here periodic ones (for $\VK_1$, there is also a right-hand side, which we do not consider for simplicity). The divergence-free constraint $\Div \Vu =0$ is necessary to guarantee the uniqueness of a solution, as otherwise the solution is only determined up to a gradient term. However, the divergence-free constraint causes some problems in the implementation of corresponding numerical methods, as e.g.\ divergence-free finite elements are difficult to find or to implement. 

With divergence-regularization, we now look for $\Vu \in \Hbf^1(Y)$ such that
\begin{equation}
\label{eq:divreg}
\int_Y \mu^{-1}\curl \Vu\cdot \curl \Vpsi^*+\Div \Vu \Div\Vpsi^*  \, dy= 0\qquad \forall \Vpsi\in \Hbf^1(Y).
\end{equation}
Clearly, any divergence-free solution to \eqref{eq:curlcurl} will also solve \eqref{eq:divreg}. On the other hand, we can insert a test function $\Vpsi= \nabla \phi$ for $\phi\in H^2(Y)$ and obtain
$$\int_Y \Div \Vu\, \triangle \phi^* \, dy =0.$$
If the domain is convex, there is $\phi\in H^2(Y)$ with $\triangle \phi = \Div \Vu$ due to elliptic regularity theory, since $\Div \Vu\in L^2(Y)$. Hence, $\Div \Vu =0$ almost everywhere.

The convexity of the domain is an essential assumption in the divergence-regularization, so that the method cannot be applied on arbitrary domains. We emphasize that divergence-regularization is only needed in the corrector here and the corresponding problems are always posed on the unit cube $Y$ or at most on a parallelepiped. Thus, convexity is guaranteed and does not impose any additional constraint. 

There are other possibilities to deal with a divergence-free constraint. The introduction of Lagrange multipliers (see \cite{CD}) leads to a mixed problem, which 
increases the computational costs and complicates the error analysis. 
The $s$-reg\-u\-lar\-i\-za\-tion suggested in \cite{DLTZ} makes the reformulation of the HMM later on (Proposition \ref{prop:HMMreform}) impossible, since different orders of derivatives appear. Thus, we choose divergence-regularization as it easily gives an equivalent problem, preserves coercivity, and can be implemented in the HMM framework as well.
\end{remark}

\begin{definition}[Cell problems and homogenized matrices]
The cell problems are to find functions $\Vv_k\in L^2(\Om;\Hbf^1_{\sharp,0}(Y))$, $v_k\in L^2(\Om; H^1_{\sharp,0}(Y))$ so that a.e.\ in $\Om$ there holds
\begin{align}
\label{eq:cellproblemmu}
&\int_Y\mu^{-1}(x, y)(\Ve_k+\curl_y\Vv_k(x,y))\cdot\curl\Vpsi^*(y)\\* \nonumber
&\qquad+\Div_y \Vv_k(x,y)\Div\Vpsi^*(y)\, dy=0 \qquad\qquad\forall\Vpsi\in \Hbf^1_{\sharp,0}(Y),\\
\label{eq:cellproblemkappa}
&\int_Y \kappa(x, y)(\Ve_k+\nabla_y v_k(x,y))\cdot \nabla \psi^*(y)\, dy=0 \quad \forall \psi\in H^1_{\sharp,0}(Y).
\end{align}
With the (unique) solutions of the cell problems \eqref{eq:cellproblemmu}--\eqref{eq:cellproblemkappa} we define the homogenized matrices
\begin{equation}
\label{eq:hommatrix}
\begin{aligned}
(\mu^{-1})^{hom}_{i,k}(x)&=\int_Y \mu^{-1}(x, y)(\de_{ik}+(\curl_y\Vv_k(x,y))_i)\, dy,\\
\kappa^{hom}_{i,k}(x)&=\int_Y \kappa(x, y)(\de_{ik}+(\nabla_y v_k(x,y))_i)\, dy, \qquad i,k=1,2,3.
\end{aligned}
\end{equation}
\end{definition}

The homogenized matrices are used to formulate the macro-scale problem for $\VE$. It has the same structure as our original problem except that the material parameters are now matrices and no scalar functions.

\begin{theorem}[Equivalence of two-scale and homogenized equation]
\label{thm:twosc-cellandmacro}
The triple $(\VE, \VK_1, K_2)$ is the unique solution of \eqref{eq:twoscE} iff $\VE\in \Hbf_0(\curl)$ solves 
\begin{equation}
\label{eq:macroE}
\int_\Om (\mu^{-1})^{hom}\curl \VE\cdot \curl \Vpsi^*-\kappa^{hom}\VE\cdot\Vpsi^*\, dx=\int_\Om \Vf\cdot\Vpsi^*\, dx 
\end{equation}
for all $\Vpsi\in \Hbf_0(\curl)$ with the matrices $(\mu^{-1})^{hom},\kappa^{hom}$ defined through \eqref{eq:hommatrix}, and with 
correctors $\VK_1, K_2$ defined as $\VK_1(x,y)=\sum_{k=1}^3(\curl\VE(x))_k\Vv_k(x, y)$,  $K_2(x,y)=\sum_{k=1}^3\VE_k(x)v_k(x, y)$, 
where $\Vv_k, v_k$ are solutions of the cell problems \eqref{eq:cellproblemmu}, \eqref{eq:cellproblemkappa}.
\end{theorem}
\begin{Proof}
Inserting the cell problems and the definition of the homogenized matrices into \eqref{eq:macroE} leads to the two-scale equation.
\end{Proof}

We end this section by a corrector-type result, which relates the two-scale solution to the asymptotic expansion. 

\begin{theorem}[Strong convergence in $\Hbf(\curl)$]
\label{th:corrector}
Let $\mu^{-1}$, $\kappa$, $\VK_1$, $\curl_x \VK_1$, $\curl_y \VK_1$, $\nabla_x K_2$ and $\nabla_y K_2$ be admissible test functions for two-scale convergence. Then it holds
$$\left\|\VE_\de(x)-\left(\VE(x)+\de\left(\VK_1\left(x, \frac{x}{\de}\right) + \nabla K_2\left(x, \frac{x}{\de}\right)\right)\right)\right\|_{\Hbf(\curl)}\TO^{\de\to0}0.$$
\end{theorem}
\begin{Proof}
Inserting the term in the norm into the heterogeneous sesquilinear form $B_\de$, using the chain rule and two-scale convergence gives the claim.
\end{Proof}
The theorem shows that the correctors $\VK_1$ and $K_2$ represent a Helmholtz decomposition of the first order term in the asymptotic expansion.  Since on the gradient subspace, the $\VH(\curl)$-norm and the $L^2$-norm are equivalent, we see that in particular $K_2$ carries important information about the solution $\VE_\de$ itself. Thus, in contrast to the elliptic case, the correctors $\VK_1$, $K_2$ have to be considered as well (and not only the weak limit $\VE$) in order to get a good $L^2$-approximation of the heterogeneous solution $\VE_\de$. This is a crucial observation. Consequently, the HMM is not only constructed to approximate $\VE$, but requires to approximate $\VK_1$ and $K_2$ as well.

\section{The Heterogeneous Multiscale Method}
\label{sec:hmm}

The basic idea of the HMM is to use a macroscopic sesquilinear form similar to \eqref{eq:macroE} for the finite element method. Instead of solving the cell problems once on the unit cube, local variants are set up and solved around the centers of the tetrahedra of some macroscopic computational grid.  In order to define the method in more detail, let us introduce some definitions.

Denote by $\CT_H=\{T_j|j\in J\}$ and $\CT_h=\{S_i|i\in I\}$ conforming, shape regular, simplicial partitions of $\Om$ and $Y$, respectively, where $\CT_h$ is additionally periodic. The $\de$-scaled and $x_j$-shifted unit cubes are denoted by $Y^\de_j=\de Y+x_j$, together with the mappings $y_j^\de:Y_j^\de\to Y$ and $x_j^\de=(y_j^\de)^{-1}: Y\to Y_j^\de$. A triangulation of the shifted unit cubes is then given by $\CT_h(Y_j^\de)=\{\tilde{S}|\tilde{S}=x_j^\de(S), S\in \CT_h\}$.
The set of interior faces is defined as 
$$\CE(\CT_H)=\{(j,l)\in J\times J|F_{jl}:=\overline{T}_j\cap\overline{T}_l\neq\emptyset, \dim(F_{jl})=2, j<l\}$$
and $\CE(\CT_h)$ with the faces $\tilde{F}_{ik}$ is defined analogously. The direct neighbors of a face $F_{jl}$ are $\hat{\om}_{F_{jl}}:=T_j\cup T_l$. The neighborhoods of vertices $V$, faces $F$, and elements $T$ are defined as
$$\om_V:=\bigcup_{j\in J, V\in \overline{T}_j}\overline{T}_j, \qquad \om_{F_{jl}}:=\bigcup_{V\in F_{jl}}\om_V, \qquad \om_{T_j}:=\bigcup_{V\in \overline{T}_j}\om_V;$$
and the neighborhood of neighbors of an element is defined as
$$\tilde{\om}_{T_j}:=\bigcup_{V\in \overline{T}_j}\bigcup_{V'\in \om_V}\om_{V'}.$$
We define the local mesh sizes $H_j:=\diam(T_j)$, $h_i:=\diam(S_i)$, $H_{jl}:=\linebreak[4]\diam(F_{jl})$, $h_{ik}:=\diam(\tilde{F}_{ik})$, and the global mesh sizes $H:=\max_{j\in J} H_j$ and $h:=\max_{i\in I}h_i$.
Finally, the discrete function spaces $\VV_{H,0}^I\subset \Hbf_0(\curl)$ and $\widetilde{W}_h^1(Y^\de_j)\subset H^1_{\sharp, 0}(Y_j^\de)$ are defined as
\begin{align*}
\VV_{H,0}^I&:=\{\Vu_H\in \Hbf_0(\curl)|\Vu_h|_T\in N_0\quad \forall T\in \CT_H\},\\
\widetilde{W}_h^1(Y_j^\de)&:=\{u_h\in H^1_{\sharp,0}(Y_j^\de)| u_h|_S\in \pz^1\quad \forall S\in \CT_h(Y_j^\de)\},
\end{align*}
where $\pz^1$ are the polynomials of maximal degree $1$ and $N_0$ is the lowest order N{\'e}d{\'e}lec element of the first family, given by $N_0:=\{\Va\times x+\Vb|\Va, \Vb\in \cz^3\}$. As in the analytical case, bold face letters indicate vector-valued functions and function spaces, for instance $\widetilde{\VW}_h^1:=(\widetilde{W}_h^1)^3$.

With these preliminaries we can now define the HMM (see also \cite{Eng,Eng1,Ohl}):
\begin{definition}[HMM]
\label{def:HMM}
The HMM-approximation of \eqref{eq:weakEde} is a discrete solution triple $(\VE_H, \VR_1(\VE_H), \VR_2(\VE_H))$, where $\VE_H\in \VV_{H, 0}^I$ is defined as the solution of
\begin{equation}
B_H(\VE_H, \Vpsi_H)=(f, \Vpsi_H) \quad \forall \Vpsi_H\in \VV_{H, 0}^I,
\end{equation}
where the discrete sesquilinear form is given by
\begin{equation}
\begin{split}
B_H(\Vu_H, \Vpsi_H):=\sum_{j\in J}\frac{|T_j|}{\de^3}&\int_{Y_j^{\de}}(\mu^{-1})_h^\de(x)\curl \VR_1(\Vu_H)(x)\cdot\curl \Vpsi_H^*(x)\\
&\quad-\kappa_h^\de(x)\VR_2(\Vu_H)(x)\cdot \VR_2(\Vpsi_H)^*(x)\, dx
\end{split}
\end{equation}
with the piecewise constant approximations $\kappa_h^\de|_{x_j^\de(S_i)}(x):=\kappa\bigl(x_j, \frac{x_j^\de(y_i)}{\de}\bigr)$ for all $S_i\in \CT_h(Y_j^\de)$ and $(\mu^{-1})_h^\de$ defined analogously. The local reconstructions $\VR_1(\Vu_H)\in \Vu_H|_{Y_j^\de}+\widetilde{\VW}_h^1(Y_j^{\de})$, $\VR_2(\Vu_H)= \Vu_H(x_k)|_{Y_j^\de}+\nabla_yu_h$ with $u_h\in \widetilde{W}_h^1(Y^{\de}_h)$ are defined as the solutions of the local cell problems
\begin{equation*}
\begin{split}
&\int_{Y_j^{\de}}(\mu^{-1})_h^\de(x)\curl \VR_1(\Vu_H)\cdot \curl \Vpsi_h^*+\Div(\VR_1(\Vu_H)-\Vu_H)\Div\Vpsi_h^*\, dx=0 \\
&\qquad \forall \Vpsi_h\in \widetilde{\VW}_h^1(Y_j^{\de})),\\
&\int_{Y_j^{\de}}\kappa_h^\de(x)\VR_2(\Vu_H)\cdot \nabla \psi_h^*\ dx=0 \quad \forall \psi_h\in \widetilde{W}_h^1(Y_j^{\de}).
\end{split}
\end{equation*}
\end{definition}

We now reformulate the reconstructions of the HMM solution triple to draw a parallel between them and the analytical correctors.

\begin{remark}[Role of the reconstructions]
\label{rem:rolecorrect}
Let $(\VE_H, \VR_1(\VE_H), \VR_2(\VE_H))$ denote the HMM-approximation from Definition \ref{def:HMM}. Setting $\VK_{j,1}(\VE_H)=\VR_{j,1}(\VE_H)-\VE_H$, we have $\VK_{j,1}\in \widetilde{\VW}_h^1(Y_j^\de)$. Furthermore, denote by $K_{j,2}(\VE_H)\in \widetilde{W}_h^1(Y_j^\de)$ the function fulfilling $\nabla K_{j,2}(\VE_H)=\VR_{j,2}(\VE_H)-\VE_H(x_j)$. We then define the discrete fine-scale corrections $\VK_{h,1}(\VE_H)\in S_H(\Om; \widetilde{\VW}_h^1(Y))$ and $K_{h,2}(\VE_H)\in S_H(\Om; \widetilde{W}_h^1(Y))$ as  
\begin{equation*}
\begin{split}
\VK_{h,1}(\VE_H)(x,y)|_{T_j\times Y}:=\frac{1}{\de}\VK_{j,1}(\VE_H)(\de y), \\
K_{h,2}(\VE_H)(x,y)|_{T_j\times Y}:=\frac{1}{\de}K_{j,2}(\VE_H)(\de y),
\end{split}
\end{equation*}
where the space of $x$-constant discrete functions is defined as
\begin{equation*}
\begin{split}
S_H(\Om; \widetilde{W}_h^1(Y))&:=\{u_h\in L^2(\Om, H^1_{\sharp,0}(Y))|u_h(\cdot, y)|_{T_j}\in \pz^0\, \forall j\in J, y\in Y\\
&\qquad \text{ and }u_h(x,\cdot)\in \widetilde{W}_h^1(Y)\, \forall x\in \Om\}.
\end{split}
\end{equation*}
The discrete fine-scale corrections $\VK_{h,1}(\VE_H)$, $K_{h,2}(\VE_H)$ are discrete counterparts of the analytical correctors $\VK_1$ and $K_2$ introduced in Theorem \ref{thm:twosc}. 
The specific relation of both will be clear from Proposition \ref{prop:HMMreform} below. Therefore, these corrections (or equivalently the reconstructions) are an important part of the HMM-approximation. As discussed at the end of Section \ref{sec:hom}, the correctors carry important information on the solution and cannot be neglected as higher order terms (in contrast to the elliptic case). In form of the fine-scale corrections, the observation transfers to the numerical scheme and the discrete setting.
\end{remark}

Having observed this correspondence, we can now reformulate the whole HMM to see that it is a direct discretization with numerical quadrature of the two-scale equation \eqref{eq:twoscE}. See \cite{Ohl} for the approach in the elliptic case.

\begin{proposition}[Reformulation of the HMM]
\label{prop:HMMreform}
Define the piecewise constant approximations $\kappa_h$ on $\Om \times Y$ by $\kappa_h(x,y)|_{T_j\times S_i}:=\kappa(x_j, y_i)$ and $\mu_h^{-1}$ in the same way. Furthermore, let $\VK_{h,1}$, $K_{h,2}$ be the discrete fine-scale corrections as defined in \ref{rem:rolecorrect}.
Then $(\VE_H, \VK_{h, 1}(\VE_H), K_{h, 2}(\VE_H))\in \VV_{H, 0}^I\times S_H(\Om; \widetilde{\VW}_h^1(Y))\times S_H(\Om; \widetilde{W}_h^1(Y))$ is a solution of
\begin{equation*}
\begin{split}
&\!\!\int_\Om \int_Y\!\!\mu^{-1}_h(x, y)(\curl\VE_H(x)+\curl_y\VK_{h,1}(\VE_H)(x, y))\!\cdot\!(\Vpsi_H^*(x)+\curl_y\Vpsi_h^*(x,y))\\
&\qquad+\Div_y\VK_{h,1}(\VE_H)(x,y)\Div_y\Vpsi_h^*(x,y)\\
&\quad\quad-\kappa(x, y)(\VE_H(x)+\nabla_y K_{h,2}(\VE_H)(x, y))\cdot(\Vpsi^*(x)+\nabla_y\psi_h^*(x,y))\, dy dx\\
&=\int_\Om \Vf(x)\cdot\Vpsi_H^*(x)\, dx\\
&\qquad\qquad\forall (\Vpsi_H, \Vpsi_h, \psi_h)\in \VV_{H, 0}^I\times L^2(\Om; \widetilde{\VW}_h^1(Y))\times L^2(\Om; \widetilde{W}_h^1(Y)).
\end{split}
\end{equation*}
\end{proposition}

\begin{Proof}
We treat the two terms of the sesquilinear form separately, but with basically the same procedure. For the first term we see from the definition of the reconstruction $\VR_1$ that for all $\Vpsi_h\in \widetilde{\VW}_h^1(Y_j^\de))$ it holds 
\begin{equation*}
\begin{split}
0&=\int_{Y_j^\de}(\mu^{-1})_h^\de(x) \curl_x (\VE_H+\VR_{j,1}(\VE_H)-\VE_H)(x)\cdot \curl_x \Vpsi_h^*(x)\\
&\qquad\quad+\Div_x(\VR_{j,1}(\VE_H)-\VE_H)(x)\Div_x \Vpsi_h^*(x)\, dx.
\end{split}
\end{equation*}
Using the transformation formula and writing $x=x^\de_j(y)$, we derive
\begin{equation*}
\begin{split}
0&=\de^3\!\int_Y \!\mu^{-1}\!\Bigl(\!x_j, \frac{x_j^\de(y)}{\de}\!\Bigr)\!\curl_x\Vpsi_h^*(x_j^\de(y))\!\cdot\!(\curl_x \VE_H(x_j)\!+\!\curl_x\VK_{j,1}(\VE_H)(x_j^\de(y)))\\
&\qquad\quad\;+\Div_x\VK_{j,1}(\VE_H)(x_j^\de(y))\Div_x\Vpsi_h^*(x_j^\de(y))\, dy,
\end{split}
\end{equation*}
since $\curl_x\VE_H(x)$ is constant on each $T_j$.
Using the definition of $\VK_{h,1}(\VE_H)$ and defining $\widetilde{\Vpsi}_h\in \widetilde{\VW}_h^1(Y)$ as $\widetilde{\Vpsi}_h(y)=\frac{1}{\de}\Vpsi_h(\de y)$, we get with the chain rule
\begin{equation*}
\begin{split}
0&=\de^3\int_Y\mu^{-1}\Bigl(x_j, \frac{x_j^\de(y)}{\de}\Bigr)\curl_y \widetilde{\Vpsi}_h^*\Bigl(\frac{x_j^\de(y)}{\de}\Bigr)\\
&\qquad\qquad\cdot \Bigl(\curl_x \VE_H(x_j)+\curl_y \VK_{h,1}(\VE_H)\Bigl(x_j, \frac{x_j^\de(y)}{\de}\Bigr)\Bigr)\\
&\qquad\quad\;+\Div_y\VK_{h,1}(\VE_H)\Bigl(x_j, \frac{x_j^\de(y)}{\de}\Bigr)\Div_y \widetilde{\Vpsi}_h^*\Bigl(\frac{x_j^\de(y)}{\de}\Bigr)\, dy.
\end{split}
\end{equation*}
As the integrand is $Y$-periodic and $\frac{x_j^\de(y)}{\de}=y+\frac{x_j}{\de}$, we finally obtain the Galerkin orthogonality
\begin{equation*}
\begin{split}
0&=\de^3\int_Y \mu^{-1}_h(x_j, y)(\curl_x \VE_H(x_j)+\curl_y\VK_{h,1}(\VE_H)(x_j, y))\cdot\curl_y\widetilde{\Vpsi}_h^*(y)\\
&\qquad\qquad+ \Div_y\VK_{h,1}(\VE_H)(x_j, y)\Div_y\widetilde{\Vpsi}_h^*(y)\, dy.
\end{split}
\end{equation*}
In the same way, we deduce
\begin{equation*}
\begin{split}
&\!\!\!\!\sum_{j\in J}\frac{|T_j|}{\de^3}\int_{Y_j^\de}(\mu^{-1})_h^\de(x)\curl_x\VR_1(\VE_H)(x)\cdot \curl \Vpsi_H^*(x)\, dx\\
&=\sum_{j\in J}|T_j| \int_Y \mu^{-1}\Bigl(x_j, \frac{x_j^\de(y)}{\de}\Bigr)\curl\Vpsi_H^*(x_j)\\
&\qquad\qquad\qquad\cdot(\curl \VE_H(x_j)+\curl_y\VK_{h,1}(\VE_H)(x_j, y))\, dy\\
&=\int_\Om\int_Y \!\mu^{-1}_h(x, y)(\curl\VE_H(x)\!+\!\curl_y\VK_{h,1}(\VE_H)(x,y))\!\cdot\!(\Vpsi_H^*(x)\!+\!\curl_y\Vpsi_h^*(x,y))\\
&\qquad\qquad\quad+\Div_y\VK_{h,1}(\VE_H)(x,y)\Div_y\Vpsi_h^*(x,y)\, dy dx.
\end{split}
\end{equation*}
In the last equality we used that the given quadrature rule is exact for the integrands and we employed the Galerkin orthogonality.
For the second term in the sesquilinear form, one can perform the same steps to reformulate the problem. For this term the computations are very similar to the elliptic case discussed in \cite[Lemma 3.5]{Ohl}.
\end{Proof}
\begin{conclusion}\label{concl:eHMM}
Let us note that the result of Theorem \ref{th:corrector} is still valid if we replace $\nabla K_2$ by $\de^{-1}\nabla_{\hspace{-2pt}y} K_2$. This implies that we can approximate $\VE_\de$ in $\Hbf(\curl)$ by 
$\VE(x)+\de \VK_1\left(x, \frac{x}{\de}\right) + \nabla_{\hspace{-2pt}y} K_2\left(x, \frac{x}{\de}\right)$.
Consequently, exploiting Proposition \ref{prop:HMMreform}, we see that our final HMM-approximation $\VE_{\mbox{\rm\tiny HMM}}$ to $\VE_\de$ can be written as
$$\VE_{\mbox{\rm\tiny HMM}}(x):=\VE_H(x)+\de \VK_{h, 1}(\VE_H)\left(x, \frac{x}{\de}\right) + \nabla_{\hspace{-2pt}y} K_{h, 2}(\VE_H)\left(x, \frac{x}{\de}\right).$$
\end{conclusion}

\section{A priori and a posteriori error analysis}
\label{sec:error}

Based on the reformulation of the HMM in Proposition \ref{prop:HMMreform}, we will give the main a priori and a posteriori error estimates in Theorems \ref{thm:aprioriH1}, \ref{thm:aprioridual}, \ref{thm:aposteriori}, and \ref{thm:lowerbound}.
All error estimates will be derived in the "two-scale energy norm"
\begin{align*}
&\|(\Vu, \Vu_1, u_2)\|_{e(G\times R)}\\*
&:=\|\curl \Vu+\curl_y \Vu_1\|_{L^2(G\times R)}+\|\Div_y \Vu_1\|_{L^2(G\times R)}+\|\Vu+\nabla_y u_2\|_{L^2(G\times R)}
\end{align*}
for $G \times R\subset \Om \times Y$ an open subdomain. If the norm is to be taken over $\Om \times Y$, we will just write $\|\cdot\|_e$. 
Let us furthermore define the error terms $e_0:=\VE-\VE_H$, $e_1:=\VK_1-\VK_{h,1}(\VE_H)$, and $e_2:=K_2-K_{h,2}(\VE_H)$.
We will only estimate these errors and leave the modeling error $\VE_\de-\left(\VE+\de\left(\VK_1\left(\cdot, \frac{\cdot}{\de}\right) + \nabla K_2\left(\cdot, \frac{\cdot}{\de}\right)\right)\right)$, introduced by homogenization, apart (cf. Theorem \ref{th:corrector} and Conclusion \ref{concl:eHMM}).

\begin{assumption}
\label{ann:Lipschitz}
On top of the periodicity of the coefficients, we also assume 
$$\mu^{-1}, \kappa\in W^{1, \infty}(\Om \times Y),$$
i.e.\ the coefficient functions are globally Lipschitz, and $\Om$ is a convex domain.  This assumption will be required for the a priori estimates (Theorems \ref{thm:aprioriH1} and \ref{thm:aprioridual}), but not for the a posteriori estimates (Theorems \ref{thm:aposteriori} and \ref{thm:lowerbound}).
\end{assumption}

\begin{theorem}[A priori estimate in the energy norm]
\label{thm:aprioriH1}
Under Assumptions \ref{assumption-coefficients} and \ref{ann:Lipschitz}, the following a priori estimate for the error between the homogenized solution and the HMM-approximation resp.\ their correctors holds:
\begin{equation*}
\|(e_0, e_1, e_2)\|_e\leq C(H+h)\|\Vf\|_{L^2(\Om)} \quad \text{ with }C=C(\Om, \kappa, \mu^{-1}).
\end{equation*}
\end{theorem}

\begin{theorem}[A priori error estimate with dual problems]
\label{thm:aprioridual}
Under the same assumptions as in Theorem \ref{thm:aprioriH1}, the Helmholtz decomposition of the error between the continuous solution $\VE$ and the HMM-approximation $\VE_H$
$$\VE-\VE_H=\nabla \theta+\Vz \quad \text{ with } \theta\in H^1_0(\Om), \Vz\bot\nabla H^1_0$$
satisfies
$$\|\theta\|_{L^2(\Om)}+\|\Vz\|_{L^2(\Om)}\leq C(H^2+h^2)\|\Vf\|_{L^2(\Om)}+C\eta_{approx}\|\Vf\|_{L^2(\Om)},$$
where $\eta_{approx}=\max\{\|\mu^{-1}-\mu^{-1}_h\|_{L^\infty(\Om\times Y)}, \|\kappa-\kappa_h\|_{L^\infty(\Om\times Y)}\}$ is a data approximation error arising from numerical quadrature.
The constant $C$ only depends on the domain $\Om$, the coefficients $\mu^{-1}$ and $\kappa$,
but not on the periodicity parameter $\de$ or the mesh sizes.
\end{theorem}

\begin{remark}
In the elliptic case, the $L^2$-norm of the error converges with quadratic rate. This better convergence is obtained by posing a dual problem and using the Aubin-Nitsche trick. The above theorem shows how the result can be transferred to problems in $\Hbf(\curl)$: On the gradient subspace, the $L^2$-norm is of the same order as the $\Hbf(\curl)$-norm, so that only on the complement a better convergence is obtained (see \cite[Remark after Thm.\ 49, p.\ 45]{Schoeberl}). Hence, the quadratic convergence here is (only) obtained in $H^{-1}$.
\end{remark}

\begin{theorem}[A posteriori error estimate]
\label{thm:aposteriori}
Let $\Vf_H$ be any piecewise polynomial approximation of $\Vf$. Under the Assumption \ref{assumption-coefficients} the error fulfills the following a posteriori error estimate
\begin{align*}
\|(e_0, e_1, e_2)\|_e&\leq C\Bigl(\sum_{j\in J}\eta_{j,1}^2+\eta_{j,2}^2\Bigr)^{1/2}+C\Bigl(\sum_{(j,l)\in \CE(\CT_H)}\eta_{jl,1}^2+\eta_{jl,2}^2\Bigr)^{1/2}\\
&\quad+C\Bigl(\sum_{j\in J}\sum_{(i,k)\in \CE(\CT_h)}\eta_{j,ik,1}^2+\eta_{j,ik,2}^2\Bigr)^{1/2}\\*
&\quad+C\Bigl(\sum_{j\in J}\zeta_j^2\Bigr)^{1/2}+C\Bigl(\sum_{j\in J}\sum_{i \in I}\zeta_{ji}^2\Bigr)^{1/2},
\end{align*}
where the constants do not depend on the mesh sizes and the periodicity parameter $\de$. 
The local indicators are defined as
\begin{align*}
\eta_{j,1}&:=H_j\Bigl\|\Vf_H+\int_Y\kappa_h(\cdot,y)(\VE_H+\nabla_y K_{h,2}(\cdot, y))\, dy\Bigr\|_{L^2(T_j)},\\
\eta_{j,2}&:=H_j\Bigl\|\Div_x\Bigl(\int_Y \kappa_h(\cdot, y)(\VE_H+\nabla_yK_{h,2}(\cdot, y))\, dy\Bigr)\Bigr\|_{L^2(T_j)},\\
\eta_{jl,1}&:=H_{jl}^{1/2}\Bigl\|\Bigl[\int_Y \mu^{-1}_h(\cdot, y)(\curl \VE_H+\curl_y\VK_{h,1}(\cdot, y))\times \Vn\, dy\Bigr]_{F_{jl}}\Bigr\|_{L^2(F_{jl})},\\
\eta_{jl,2}&:=H_{jl}^{1/2}\Bigl\|\Bigl[\int_Y \kappa_h(\cdot, y)(\VE_H+\nabla_y K_{h,2}(\cdot, y))\cdot \Vn\, dy\Bigr]_{F_{jl}}\Bigr\|_{L^2(F_{jl})}\\
\eta_{j,ik,1}&:=h_{ik}^{1/2}\|[\mu^{-1}_h(\curl \VE_H+\curl_y \VK_{h,1})\times \Vn+\Div_y \VK_{h,1}\Vn]_{\tilde{F}_{ik}}\|_{L^2(T_j\times \tilde{F}_{ik})},\\
\eta_{j,ik,2}&:=h_{ik}^{1/2}\|[\kappa_h(\VE_H+\nabla_y K_{h,2})\cdot \Vn]_{\tilde{F}_{ik}}\|_{L^2(T_j\times\tilde{F}_{ik})},\\
\zeta_j&:=H_j\|\Vf_H-\Vf\|_{L^2(T_j)},\\
\zeta_{ji}&:=\|(\mu^{-1}_h-\mu^{-1})(\curl\VE_H+\curl_y \VK_{h,1})\|_{L^2(T_j\times S_i)}\\
&\quad+\|(\kappa_h-\kappa)(\VE_H+\nabla_y K_{h,2})\|_{L^2(T_j\times S_i)}.
\end{align*}
Here and in the following, $[\cdot]_F$ denotes the jump across the face $F$.
\end{theorem}

\begin{remark}[Discussion of the error indicators]
The error indicators can be split into two groups: $\zeta_j$ and $\zeta_{ji}$ are data approximation errors, which come from the use of numerical quadrature. The error indicators denoted by $\eta$ are different contributions to the discretization error: $\eta_{j,1}$ is the element residual on the macro-scale, $\eta_{jl}$ are the jump residuals on the macro-scale (in normal and tangential direction), and $\eta_{j,ik}$  are the jump residuals on the micro-scale. $\eta_{j,2}$ indicates how well the deduced equation $\Div_x (\kappa\VE_0)=0$ is fulfilled in the discrete case. Here, the assumption $\Div \Vf =0$ (in the weak sense) has an effect on the error estimator: If we just assume $\Vf\in \VH(\Div)$, the deduced equation is $\Div_x(\kappa \VE_0+\Vf)=0$ and thus we have additional terms $\Div_x \Vf_H$ in $\eta_{j,2}$ and $\Vf_H\cdot \Vn$ in $\eta_{jl,2}$ for the polynomial approximation $\Vf_H$ of $\Vf$. Furthermore, in the data approximation error $\zeta_j$ we then have to take the $\VH(\Div, T_j)$-norm.
\end{remark}

\begin{theorem}[Lower bound on the error]
\label{thm:lowerbound}
With the same notations and under the same assumptions as in Theorem \ref{thm:aposteriori}, the following local bounds on the error hold:
\begin{align*}
\eta_{j,1}&\leq C\Bigl(\|(e_0, e_1, e_2)\|_{e(T_j\times Y)}+\zeta_j+\Bigl(\sum_{i}\zeta_{ji}^2\Bigr)^{1/2}\Bigr),\\
\eta_{j,2}&\leq C\Bigl(\|(e_0, e_1, e_2)\|_{e(T_j\times Y)}+\Bigl(\sum_{i}\zeta_{ji}^2\Bigr)^{1/2}\Bigr),\\
\eta_{jl,1}&\leq C\Bigl(\|(e_0, e_1, e_2)\|_{e(\hat{\om}_{F_{jl}}\times Y)}+(\zeta_j^2+\zeta_l^2)^{1/2}+\Bigl(\sum_i\zeta_{ji}^2+\zeta_{li}^2\Bigr)^{1/2}\Bigr),\\
\eta_{jl,2}&\leq C\Bigl(\|(e_0, e_1, e_2)\|_{e(\hat{\om}_{F_{jl}}\times Y)}+\Bigl(\sum_i\zeta_{ji}^2+\zeta_{li}^2\Bigr)^{1/2}\Bigr),\\
\eta_{j,ik,\nu}&\leq C(\|(e_0, e_1, e_2)\|_{e(T_j\times \hat{\om}_{\tilde{F}_{ik}})}+(\zeta_{ji}^2+\zeta_{jk}^2)^{1/2}), \quad \nu=1,2.
\end{align*}
Thus, we have the global estimate
\begin{equation*}
\begin{split}
&\!\!\!\!\!\!\!\!\Bigl(\sum_{j\in J}\eta_{j,1}^2+\eta_{j,2}^2\Bigr)^{1/2}+\Bigl(\sum_{(j,l)\in \CE(\CT_H)}\eta_{jl,1}^2+\eta_{jl,2}^2\Bigr)^{1/2}\\
&\quad+\Bigl(\sum_{j\in J}\sum_{(i,k)\in \CE(\CT_h)}\eta_{j,ik,1}^2+\eta_{j,ik,2}^2\Bigr)^{1/2}\\
&\leq C\Bigl(\|(e_0, e_1, e_2)\|_e+\Bigl(\sum_{j \in J}\zeta_j^2\Bigr)^{1/2}+\Bigl(\sum_{j\in J}\sum_{i\in I}\zeta_{ji}^2\Bigr)^{1/2}\Bigr).
\end{split}
\end{equation*}
\end{theorem}

\begin{remark}
Theorems \ref{thm:aposteriori} and \ref{thm:lowerbound} together show that the local error indicators are reliable and efficient with respect to the two-scale homogenized problem from Theorem \ref{thm:twosc}. Up to data approximation errors, the error and the indicators converge with the same rate. Thus, the indicators can be used for adaptive algorithms, e.g.\ for mesh refinement, both on the coarse and fine scale (cf. \cite{Ohl} for related adaptive algorithms in the elliptic case). 
\end{remark}

\section{Proofs of the main results}
\label{sec:errorproof}
In this section, the essential proofs of the homogenization result (namely the two-scale equation) and the error estimates for the HMM will be given.

\begin{Proof}[Proof of Theorem \ref{thm:twosc}]
As solutions to \eqref{eq:weakEde} are uniformly bounded in $\Hbf(\curl)$, we have by Theorem \ref{thm:Visintin} that, along a subsequence, $\VE_\de\twosc\VE_0$ in $(L^2(\Om\times Y))^3$ and $\curl\VE_\de\twosc\curl\VE+\curl_y\VK_1$. Furthermore, we have from Theorem \ref{thm:Well} that $\VE_0=\VE+\nabla_yK_2$, where $\VE=\int_Y\VE_0(\cdot,y)\, dy$ is the weak limit of $\VE_\de$ in $\Hbf(\curl)$. We insert $\Vpsi(x)=\Vw(x)+\de\Vw_1(x,\frac{x}{\de})$ with arbitrary $\Vw \in (C^{\infty}_0(\Om))^3, \Vw_1\in C^{\infty}_0(\Om; (C^{\infty}_\sharp(Y))^3)$ as a test function in \eqref{eq:weakEde}. 
Because of the assumption \eqref{eq:assptstrong2sc} on the parameters we can apply two-scale convergence to each of the terms in \eqref{eq:weakEde} and thereby obtain
\begin{equation}
\label{eq:twosc1}
\begin{split}
&\!\!\int_\Om\int_Y\mu^{-1}(x, y)(\curl\VE(x)+\curl_y\VK_1(x,y))\cdot(\curl\Vw^*(x)+\curl_y\Vw_1^*(x,y))\\
&\qquad \quad -\kappa(x, y)(\VE(x)+\nabla_yK_2(x,y))\cdot\Vw^*(x)\, dy dx\\
&=\int_\Om \Vf(x)\cdot\Vw^*(x)\, dx.
\end{split}
\end{equation}
By density this also holds for test functions in $\Hbf_0(\curl)\times L^2(\Om; \Hbf^1_{\sharp,0}(Y))$.
Furthermore, we can deduce from \eqref{eq:weakEde} with $\Div \Vf=0$ that it holds
$$\de\int_\Om\kappa_\de(x)\VE_\de(x)\cdot \nabla \psi^*(x)\, dx=0 \quad \forall\psi\in H^1_0(\Om).$$
Choosing $\psi(x, \frac{x}{\de})\in C^\infty_0(\Om; C^{\infty}_\sharp(Y))$, we obtain with two-scale convergence 
$$ \int_\Om\int_Y \kappa(x, y)(\VE(x)+\nabla_yK_2(x,y))\cdot \nabla_y\psi^*(x,y)\, dy dx=0.$$
Inserting this into \eqref{eq:twosc1}, we get \eqref{eq:twoscE} except for the divergence term, but with the additional constraint $\Div_y \VK_1=0$.

As discussed in Remark \ref{rem:divreg}, we can apply divergence-regularization in this case to obtain an equivalent problem. Looking at the method, we directly see that \eqref{eq:twoscE} is simply the regularization of \eqref{eq:twosc1}. The equivalence of the problems can be seen as discussed in the remark, just insert $\Vpsi=0$, $\psi_2=0$, and $\Vpsi_1 =\nabla \varphi$ with $\varphi \in L^2(\Om; H^2_{\sharp,0}(Y))$ as test function in \eqref{eq:twosc1}.

So far we have shown \eqref{eq:twoscE} just for a subsequence. If we can prove that the solution of the two-scale equation is unique, the result holds for the whole sequence.
In fact we will prove that \eqref{eq:twoscE} is of the form $\CB(u, \psi)=\VF(\psi)$ with a continuous and coercive sesquilinear form $\CB$ and a functional $\VF$. As for the weak solutions to \eqref{eq:weakEde}, Lax-Milgram-Babu{\v{s}}ka then yields the uniqueness of the solution. Moreover, this reformulation will be important for the error estimates later on.
We consider the Hilbert space 
$$\Hbf:=\Hbf_0(\curl)\times L^2(\Om; \Hbf^1_{\sharp,0}(Y))\times L^2(\Om; H^1_{\sharp,0}(Y))$$
with its natural norm 
$$ \|(\Vu, \Vu_1, u_2)\|_\VH^2=\|\Vu\|_{\VH(\curl)}^2+\|\Vu_1\|^2_{L^2(\Om; \Hbf^1(Y))}+\|u_2\|^2_{L^2(\Om; H^1(Y))}.$$
Clearly, the right-hand side is in the dual space of $\Hbf$ and the left-hand side defines a continuous sesquilinear form $\CB$. With the same computations as for the existence of  a weak solution, one can also show that $\CB$ is coercive with respect to the energy norm
$$\|\curl\Vu+\curl_y\Vu_1\|^2_{L^2(\Om\times Y)}+\|\Div_y \Vu_1\|^2_{L^2(\Om\times Y)}+\|\Vu+\nabla_yu_2\|^2_{L^2(\Om\times Y)}.$$
It remains to show the equivalence of the energy and the natural norm.
It holds
\begin{equation*}
\begin{split}
&\!\!\!\!\|\curl \Vu+\curl_y\Vu_1\|^2_{L^2(\Om\times Y)}\\
&=\int_\Om\int_Y |\curl\Vu|^2+|\curl_y\Vu_1|^2+2\Re(\curl_y\Vu_1(x,y)\cdot \curl\Vu^*(x))\, dy dx\\
&=\|\curl\Vu\|^2_{L^2(\Om)}+\|\curl_y\Vu_1\|^2_{L^2(\Om\times Y)}\\
&\quad-2\Re\underbrace{\left(\int_\Om\int_{\partial Y}(\Vu_1(x,y)\times \Vn)\cdot \Vu^*(x)\, d\sigma dx\right)}_{=0, \text{ periodicity of }\Vu_1}\\
&=\|\curl\Vu\|^2_{L^2(\Om)}+\|\curl\Vu_1\|^2_{L^2(\Om\times Y)}.
\end{split}
\end{equation*}
With integration by parts and the Poincar{\'e} inequality we see that $\|\curl_y\Vu_1\|_{L^2}+\|\Div_y\Vu_1\|_{L^2}$ is equivalent to the full $H^1$-norm. 
Similarly, we derive
\begin{equation*}
\begin{split}
&\!\!\!\!\|\Vu+\nabla_yu_2\|^2_{L^2(\Om\times Y)}\\
&=\|\Vu\|^2_{L^2(\Om)}+\|\nabla_yu_2\|^2_{L^2(\Om\times Y)}+2\Re\underbrace{\left(\int_\Om \int_{\partial Y}\Vu(x)\cdot \Vn \,u_2^*(x,y)\, d\sigma dx\right)}_{=0, \text{ periodicity of }u_2},
\end{split}
\end{equation*}
and again with the Poincar{\'e} inequality the $H^1$-seminorm is equivalent to the full $H^1$-norm. Due to the uniqueness of the two-scale solution the whole sequence $\VE_\de$ converges as asserted in the theorem.  
\end{Proof}

Having identified the variational formulations in Theorem \ref{thm:twosc} and Proposition \ref{prop:HMMreform}, we collect some useful notation for the proofs of the error estimates in the next paragraphs.

\begin{notation}[Sesquilinear forms $\CB$, $\CB_h$ and residual]
We define the continuous and discrete sesquilinear forms $\CB, \CB_h: [\Hbf_0(\curl)\times L^2(\Om; \Hbf^1_{\sharp,0}(Y) )\times L^2(\Om; H^1_{\sharp,0}(Y))]^2\to \cz$ as 
\begin{align*}
&\!\!\!\!\CB((\VE, \VK_1, K_2), (\Vpsi, \Vpsi_1, \psi_2)):=\\*
&\int_\Om\int_Y \mu^{-1}(x,y)(\curl\VE(x)+\curl_y\VK_1(x,y))\cdot(\curl\Vpsi^*(x)+\curl_y\Vpsi_1^*(x,y))\\
&\qquad\quad+\Div_y \VK_1(x,y)\Div_y \Vpsi_1^*(x,y)\\*
&\qquad\quad-\kappa(x,y)(\VE(x)+\nabla_y K_2(x,y))\cdot(\Vpsi^*(x)+\nabla_y \psi_2^*(x,y))\, dy dx,\\
&\!\!\!\!\CB_h((\Vu_H, \Vu_h, u_h), (\Vpsi_H, \Vpsi_h, \psi_h)):=\\*
&\int_\Om\int_Y \mu^{-1}_h(x,y)(\curl\Vu_H(x)+\curl_y\Vu_h(x,y))\cdot(\curl\Vpsi_H^*(x)+\curl_y\Vpsi_h^*(x,y))\\
&\qquad\quad+\Div_y \Vu_h(x,y)\Div_y \Vpsi_h^*(x,y)\\*
&\qquad\quad-\kappa_h(x,y)(\Vu_H(x)+\nabla_y u_h(x,y))\cdot(\Vpsi_H^*(x)+\nabla_y \psi_h^*(x,y))\, dy dx.
\end{align*}
In addition we define the residual 
\begin{align}
\nonumber
\Res_h:\Hbf_0(\curl, \Om)\times L^2(\Om; \Hbf^1_{\sharp,0}(Y))&\times L^2(\Om; H^1_{\sharp,0}(Y))\\\nonumber
&\hspace{-4cm}\to \Hbf_0(\curl)'\times L^2(\Om; (\VH^{-1}_{\sharp,0}(Y))^3)\times L^2(\Om; H^{-1}_{\sharp,0}(Y)) \\
\label{eq:residual}
\text{ as }\quad\langle \Res_h(\Vu, \Vu_1, u_2),(\Vpsi, \Vpsi_1, \psi_2)\rangle&:=\CB_h((\Vu, \Vu_1, u_2), (\Vpsi, \Vpsi_1, \psi_2))-(\Vf, \Vpsi).
\end{align}
\end{notation}

In the following, we will write $\VK_{h,1}$ instead of $\VK_{h,1}(\VE_H)$, and $K_{h,2}$ instead of $K_{h,2}(\VE_H)$. If it is clear on which variables ($x$, $y$) functions depend, we will omit those variables for the sake of readability. $C$ denotes a generic constant, independent of the mesh sizes and $\de$.

\subsection{Proofs of the a priori estimates}

The a priori estimates are based on the C{\'e}a lemma, dual problems and interpolation operators. As the assumptions on the coefficients and the domain imply higher (namely $H^2$) regularity of the solution, we can use the nodal interpolation operators.

\begin{lemma}[Lagrange interpolation operator]
\label{lem:lagrangeintpol}
Denote by $I_h^L: C^0(\overline{Y})\to W_h^1(Y)$ the standard Lagrange interpolation operator. Define now $\tilde{I}_h^L$ for $v\in L^2(\Om; C^0(\overline{Y}))$ by
$$\tilde{I}_h^L(v)(x,y):=(I_h^L(v)(x, \cdot))(y)-\int_Y (I_h^L(v)(x, \cdot))(s)\, ds.$$
For all $v \in L^2(\Om; C^0(\overline{Y}))\cap L^2(\Om; H^1_{\sharp,0}(Y))$ this interpolation operator is well-defined with $\tilde{I}_h^L(v)\in L^2(\Om; \widetilde{W}_h^1(Y))$ and the estimate
\begin{equation*}
\|\tilde{I}_h^L(v)-v\|_{L^2(\Om\times Y)} + h\, |\tilde{I}_h^L(v)-v|_{L^2(\Om, H^1(Y))}\leq Ch^2|v|_{L^2(\Om, H^2(Y))}.
\end{equation*}
For vector functions the interpolation operator is defined component-wise and fulfills the same estimates.
\end{lemma}
\begin{Proof}
For a proof we refer the reader to \cite{Henn}.
\end{Proof}

\begin{lemma}[Edge interpolation operator]
\label{lem:edgeintpol}
Denote by $I_H^E:\Hbf^1(\curl)\to \VV_{H,0}^I$ the nodal interpolation operator for the N{\'e}d{\'e}lec elements. It fulfills the estimate
\begin{equation*}
\|\Vu-I_H^E(\Vu)\|_{\Hbf(\curl, T_j)}\leq CH_j\|\Vu\|_{\Hbf^1(\curl, T_j)}.
\end{equation*}
\end{lemma}
\begin{Proof}
For the definition of $I_H^E$ and the proof we refer to \cite{Monk}.
\end{Proof}

\begin{Proof}[Proof of Theorem \ref{thm:aprioriH1}]
Denote by $\widetilde{\VE}_H\in \VV_{H, 0}^I$, $\widetilde{\VE}_h\in L^2(\Om; \widetilde{\VW}_h^1(Y))$, and $\widetilde{E}_h\in L^2(\Om; \widetilde{W}_h^1(Y))$ the unique solution of 
\begin{equation*}
\begin{split}
&\CB((\widetilde{\VE}_H, \widetilde{\VE}_h, \widetilde{E}_h), (\Vpsi_H, \Vpsi_h, \psi_h))=(\Vf, \Vpsi_H)\\
&\qquad\forall \Vpsi_H\in \VV_{H, 0}^I, \Vpsi_h\in L^2(\Om; \widetilde{\VW}_h^1(Y)), \psi_h\in L^2(\Om; \widetilde{W}_h^1(Y)).
\end{split}
\end{equation*}
According to C{\'e}a's Lemma it holds
\begin{equation*}
\begin{split}
&\|(\VE-\widetilde{\VE}_H, \VK_1-\widetilde{\VE}_h, K_2-\widetilde{E}_h)\|_e\\
&\leq C\bigl(\inf_{\Vpsi_H \in \VV_{H, 0}^I}\|\VE-\Vpsi_H\|_{\Hbf(\curl)}+\inf_{\Vpsi_h\in L^2(\Om; \widetilde{\VW}_h^1(Y))}|\VK_1-\Vpsi_h|_{L^2(\Om; \Hbf^1(Y))}\\
&\qquad\quad+\inf_{\psi_h\in L^2(\Om; \widetilde{W}_h^1(Y))}|K_2-\psi_h|_{L^2(\Om; H^1(Y))}\bigr).
\end{split}
\end{equation*}
With the interpolation estimates of Lemmas \ref{lem:lagrangeintpol} and \ref{lem:edgeintpol} we now derive
\begin{equation*}
\begin{split}
&\!\!\!\!\|(\VE-\widetilde{\VE}_H, \VK_1-\widetilde{\VE}_h, K_2-\widetilde{E}_h)\|_e\\
&\leq C(\|\VE-I_H^E(\VE)\|_{\Hbf(\curl)}+|\VK_1-\tilde{I}_h^L(\VK_1)|_{L^2(\Om; \Hbf^1(Y))}\\
&\qquad\;+|K_2-\tilde{I}_h^L(K_2)|_{L^2(\Om; H^1(Y))})\\
&\leq C(H\|\VE\|_{\Hbf(\curl)}+h|\VK_1|_{L^2(\Om; \VH^2(Y))}+h|K_2|_{L^2(\Om; H^2(Y))})\\
&\leq C(H+h)\|\Vf\|_{L^2(\Om)},
\end{split}
\end{equation*}
where in the last inequality we used regularity and stability results for the analytic solution. (Note that because of our assumptions on the parameters and on $\Om$ the two-scale solution admits $H^2$ regularity.) 
Furthermore, because of the definition of $(\widetilde{\VE}_H, \widetilde{\VE}_h, \widetilde{E}_h)$ it holds
\begin{equation*}
\begin{split}
&\!\!\!\!\|(\widetilde{\VE}_H-\VE_H, \widetilde{\VE}_h-\VK_{h,1}(\VE_H), \widetilde{E}_h-K_{h,2}(\VE_H))\|^2_e\\
&\leq C \bigl|(\CB_h-\CB)((\VE_H, \VK_{h,1}(\VE_H), K_{h,2}(\VE_H)),\\
&\hspace{2.4cm}(\widetilde{\VE}_H-\VE_H, \widetilde{\VE}_h-\VK_{h,1}(\VE_H), \widetilde{E}_h-K_{h,2}(\VE_H)))\bigr|\\
&\leq C\max\{\|\mu_h^{-1}-\mu^{-1}\|_{L^\infty(\Om \times Y)}, \|\kappa_h-\kappa\|_{L^\infty(\Om \times Y)}\}\\
&\qquad\|(\VE_H, \VK_{h,1}(\VE_H), K_{h,2}(\VE_H))\|_e\\
&\qquad\|(\widetilde{\VE}_H-\VE_H, \widetilde{\VE}_h-\VK_{h,1}(\VE_H), \widetilde{E}_h-K_{h,2}(\VE_H))\|_e.
\end{split}
\end{equation*} 
From the Lipschitz continuity (with constant $L$) it follows
$$\|\kappa_h-\kappa\|_{L^\infty(\Om\times Y)}\leq L \sup_{(x,y)\in \Om \times Y}|(x_i, y_k)-(x,y)|\leq L(H+h),$$
and the same estimate also applies to $\|\mu^{-1}_h-\mu^{-1}\|_{L^\infty(\Om\times Y)}$.
Together with a stability estimate for the HMM approximation this yields 
$$\|(\widetilde{\VE}_H-\VE_H, \widetilde{\VE}_h-\VK_{h,1}(\VE_H), \widetilde{E}_h-K_{h,2}(\VE_H))\|_e\leq C(H+h)\|\Vf\|_{L^2(\Om)}.$$
Splitting the total error $\VE-\VE_H$ into the contributions $\VE-\widetilde{\VE}_H$ and $\widetilde{\VE}_H-\VE_H$ and using the two estimates, we obtain the assertion.
\end{Proof}

\begin{Proof}[Proof of Theorem \ref{thm:aprioridual}]
As the terms in the Helmholtz decomposition are orthogonal w.r.t.\ the $L^2$-scalar product, we have 
\begin{align}
\label{helmholtz-decomp-est}
\|\nabla\theta\|_{L^2}\leq\|\nabla \theta\|_{L^2}+\|\Vz\|_{L^2}= \|\VE-\VE_H\|_{L^2}.
\end{align}

To estimate $\Vz$, let $(\Vw, \Vw_1, w_2)$ be the solution of the dual problem
$$\CB((\Vpsi, \Vpsi_1, \psi_2), (\Vw, \Vw_1, w_2))=(\Vz, \Vpsi)_{L^2} \quad \forall (\Vpsi, \Vpsi_1, \psi_2),$$
and $(\Vw_H, \Vw_h, w_h)$ the solution of the corresponding discrete dual problem
$$\CB_h((\Vpsi_H, \Vpsi_h, \psi_h), (\Vw_H, \Vw_h, w_h))=(\Vz, \Vpsi_H)_{L^2} \quad \forall (\Vpsi_H, \Vpsi_h, \psi_h).$$
The analytical and discrete spaces are the same as in the problems for $\VE$ and $\VE_H$ and therefore not given again here.
Because of $(\nabla\theta, \Vz)=0$ it holds $\|\Vz\|^2_{L^2}=(\Vz, \VE-\VE_H)$. Thus, it follows
$$\|\Vz\|^2_{L^2}=(\Vz, \VE-\VE_H)_{L^2}=\CB((e_0,e_1,e_2), (\Vw, \Vw_1, w_2)).$$
Using the definition of $\VE$ as exact solution and of $\VE_H$ as the HMM-approx\-imation, we deduce
\begin{equation*}
\begin{split}
\|\Vz\|_{L^2}^2&=\CB((e_0,e_1,e_2), (\Vw, \Vw_1, w_2))-\CB((\VE, \VK_1, K_2), (\Vw_H, \Vw_h, w_h))\\
&\quad+\CB_h((\VE_H, \VK_{h,1}(\VE_H), K_{h,2}(\VE_H)), (\Vw_H, \Vw_h, w_h))\\
&=\CB((e_0, e_1, e_2),(\Vw-\Vw_H, \Vw_1-\Vw_h, w_2-w_h))\\
&\quad+(\CB_h-\CB)((\VE_H, \VK_{h,1}(\VE_H), K_{h,2}(\VE_H)), (\Vw_H, \Vw_h, w_h))\\
&\leq C\|(e_0, e_1, e_2)\|_e\|(\Vw-\Vw_H, \Vw_1-\Vw_h, w_2-w_h)\|_e\\
&\quad+C\max\{\|\mu^{-1}-\mu^{-1}_h\|_{L^\infty(\Om\times Y)}, \|\kappa-\kappa_h\|_{L^\infty(\Om\times Y)}\}\\
&\qquad\qquad\|(\VE_H, \VK_{h,1}(\VE_H), K_{h,2}(\VE_H))\|_e\|(\Vw_H, \Vw_h, w_h)\|_e.
\end{split}
\end{equation*}
According to Theorem \ref{thm:aprioriH1}, it holds 
$$\|(\Vw-\Vw_H, \Vw_1-\Vw_h, w_2-w_h)\|_e\leq C (H+h)\|\Vz\|_{L^2}.$$
Hence, together with stability estimates for $\VE_H$ and $\Vw_H$ it follows
$$\|\Vz\|_{L^2}^2\leq C(H+h)\|\Vz\|_{L^2}\|(e_0, e_1, e_2)\|_e+C\eta_{approx}\|\Vz\|_{L^2} \| \Vf\|_{L^2}.$$

To estimate $\theta$, we pose another dual problem: Find $\hat{w}\in H^1_0(\Om)$, $\hat{w}_2\in L^2(\Om; H^1_{\sharp,0}(Y))$ such that
$$\CA((\hat{\psi}, \hat{\psi}_2), (\hat{w}, \hat{w}_2))=(\theta, \hat{\psi})_{L^2} \quad \forall (\hat{\psi}, \hat{\psi}_2)\in H^1_0(\Om)\times L^2(\Om; H^1_{\sharp,0}(Y))$$
with
\begin{equation*}
\begin{split}
\CA(&(\hat{\psi}, \hat{\psi}_2), (\hat{u}, \hat{u}_2))\\
&:=-\int_\Om \int_Y \kappa(x,y)(\nabla \hat{\psi}(x)+\nabla_y \hat{\psi}_2(x,y))\cdot(\nabla \hat{u}^*(x)+\nabla_y \hat{u}^*_2(x,y))\, dy dx.
\end{split}
\end{equation*}
Again, let us denote by $(\hat{w}_H, \hat{w}_h)\in W_H^1(\Om)\times L^2(\Om; \widetilde{W}_h^1(Y))$ the solution of the corresponding discrete dual problem. This dual problem is related to our original problem by the equation
$$\CA((\hat{\psi}, \hat{\psi}_2), (\hat{w}, \hat{w}_2))=\CB((\nabla \hat{\psi}, \Vpsi_1, \hat{\psi}_2), (\nabla \hat{w}, 0, \hat{w}_2))$$
for all $(\hat{\psi}, \Vpsi_1, \hat{\psi}_2)\in H^1_0(\Om)\times L^2(\Om; \Hbf^1_{\sharp,0}(Y))\times L^2(\Om; H^1_{\sharp,0}(Y))$.
Inserting $\hat{\psi}=\theta$ and $\hat{\psi}_2=e_2$, we then obtain
\begin{equation*}
\begin{split}
\|\theta\|_{L^2}^2&=\CA((\theta, e_2), (\hat{w}, \hat{w}_2))=\CB((\nabla \theta, e_1, e_2), (\nabla \hat{w}, 0, \hat{w}_2))\\
&=\CB((e_0, e_1, e_2), (\nabla \hat{w}, 0, \hat{w}_2))-\CB((\Vz, 0, 0), (\nabla \hat{w}, 0, \hat{w}_2)).
\end{split}
\end{equation*}
By the properties of the Helmholtz decomposition we have $(\Vz, \nabla \hat{w}_H)_{L^2}=0$.
With the same computations as for the dual problem with $\Vz$, we then derive
\begin{equation*}
\begin{split}
\|\theta\|_{L^2}^2&=\CB((\overset{=\nabla \theta}{\overbrace{e_0-\Vz}}, e_1, e_2), (\nabla(\hat{w}-\hat{w}_H), 0, \hat{w}_2-\hat{w}_h))\\
&\quad-\CB((\Vz, 0, 0), (\nabla \hat{w}, 0, \hat{w}_2))\\
&\quad+(\CB_h-\CB)(\VE_H, \VK_{h,1}(\VE_H), K_{h,2}(\VE_H)), (\nabla \hat{w}_H, 0, \hat{w}_h))\\
&\leq \bigl|\CA((\theta, e_2)(\hat{w}-\hat{w}_H, \hat{w}_2-\hat{w}_h))\bigr|+C\|\Vz\|_{L^2}\|\nabla \hat{w}+\nabla_y \hat{w}_2\|_{L^2(\Om\times Y)}\\
&\quad+C\max\{\|\mu^{-1}-\mu^{-1}_h\|_{L^\infty(\Om\times Y)}, \|\kappa-\kappa_h\|_{L^\infty(\Om\times Y)}\}\\
&\qquad\quad\|(\VE_H, \VK_{h,1}(\VE_H), K_{h,2}(\VE_H))\|_e\|\nabla \hat{w}_H+\nabla_y \hat{w}_h\|_{L^2(\Om\times Y)}\\
&\leq C\|\nabla \theta+\nabla_y e_2\|_{L^2(\Om\times Y)}\|\nabla (\hat{w}-\hat{w}_H)+\nabla_y (\hat{w}_2-\hat{w}_h)\|_{L^2(\Om\times Y)}\\
&\quad+C\eta_{approx}\|\theta\|_{L^2} \| \Vf\|_{L^2}+C\|\Vz\|_{L^2}\|\theta\|_{L^2},
\end{split}
\end{equation*}
where in the last inequality we used the stability estimate for $\VE_H$ and a stability estimate for the solution of elliptic two-scale equations.
From a priori error estimates for elliptic two-scale problems (see \cite{Henn, Ohl}), we know that 
$$\|\nabla(\hat{w}-\hat{w}_H)+\nabla_y (\hat{w}_2-\hat{w}_h)\|\leq C(H+h)\|\theta\|_{L^2}.$$
Inserting this, the estimate for the Helmholtz decomposition \eqref{helmholtz-decomp-est}, and the estimate for $\Vz$ from above, we finally obtain
\begin{equation*}
\begin{split}
\|\theta\|_{L^2}^2&\leq C(H+h)\|\theta\|_{L^2}\underbrace{(
\|e_0\|_{L^2(\Omega)}
+\|\nabla_y e_2\|_{L^2}}_{\leq \|(e_0, e_1, e_2)\|_e})\\
&\quad+C(H+h)\|(e_0, e_1, e_2)\|_e\|\theta\|_{L^2}+C\eta_{approx}\|\theta\|_{L^2} \| \Vf\|_{L^2}.
\end{split}
\end{equation*}

The estimates for $\Vz$ and $\theta$ together with the a priori error estimate of Theorem \ref{thm:aprioriH1} give the claim.
\end{Proof}

\subsection{Proofs of the a posteriori estimates}

For the a posteriori estimates we no longer assume higher regularity and therefore, need other interpolation operators. 

\begin{lemma}[Error estimates for the Cl{\'e}ment interpolation operator]
\label{lem:clementintpol} 
Let us denote by $\bar{I}_h:L^2(Y)\to \widetilde{W}_h^1(Y)$ the Cl{\'e}ment interpolation operator (see \cite{C}) with appropriate adaptations to periodic boundary conditions and zero average. We define $I_h:L^2(\Om; L^2(Y))\to L^2(\Om; \widetilde{W}_h^1(Y))$ as
$$I_hu(x,y):=\overline{I}_h(u(x, \cdot))(y)\quad \forall x \in \Om.$$ 
Then the following estimates hold for all $u\in L^2(\Om; H^1_{\sharp,0}(Y))$:
\begin{equation*}
\begin{split}
\|u-I_hu\|_{L^2(T_j\times S_i)}&\leq C h_i \|\nabla_y u\|_{L^2(T\times \om_{S_i})},\\
\|u-I_hu\|_{L^2(T_j\times \tilde{F}_{ik})}&\leq C h^{1/2}_{ik} \|\nabla_y u\|_{L^2(T\times \om_{\tilde{F}_{ik}})}.
\end{split}
\end{equation*}
Again, the Cl{\'e}ment operator can be defined component-wise for vector functions and then fulfills the same estimates.
\end{lemma}
\begin{Proof}
A proof can be found in \cite{C}.
\end{Proof}

\begin{lemma}[Sch{\"o}berl interpolation operator]
\label{lem:schoeberlintpol}
There exists an operator $I_H: \Hbf_0(\curl)\to \VV_{H,0}^I$ such that for every $\Vu\in \Hbf_0(\curl)$ there exist $\theta\in H^1_0(\Om)$ and $\Vz\in \VH^1_0(\Om)$ with 
$$\Vu-I_H\Vu=\nabla \theta+\Vz.$$
The decomposition fulfills the estimates 
\begin{equation*}
\begin{split}
H_j^{-1}\|\theta\|_{L^2(T_j)}+\|\nabla\theta\|_{L^2(T_j)}&\leq C\|\Vu\|_{L^2(\tilde{\om}_{T_j})},\\
H_j^{-1}\|\Vz\|_{L^2(T_j)}+\|\nabla \Vz\|_{L^2(T_j)}&\leq C \|\curl \Vu\|_{L^2(\tilde{\om}_{T_j})}.
\end{split}
\end{equation*}
Together with the trace inequality we moreover have the estimates
\begin{equation*}
\begin{split}
\|\theta\|_{L^2(F_{jl})}&\leq C H_{jl}^{1/2}\|\Vu\|_{L^2(\tilde{\om}_{T_j})},\\
\|\Vz\|_{L^2(F_{jl})}&\leq C H_{jl}^{1/2}\|\curl\Vu\|_{L^2(\tilde{\om}_{T_j})},
\end{split}
\end{equation*}
where, of course, $T_j$ can also be substituted by $T_l$.
\end{lemma}
\begin{Proof}
For the construction of $I_H$ and a proof of the estimates we refer to \cite{Sch1}. Additional details on $I_H$ can also be found in \cite{Sch2, Sch3}.
\end{Proof}

\begin{Proof}[Proof of Theorem \ref{thm:aposteriori}]

First of all, we derive an error identity.
From the definition of the error terms and Proposition \ref{prop:HMMreform} we deduce
\begin{equation*}
\begin{split}
&\!\!\!\!\CB((e_0, e_1, e_2), (\Vpsi, \Vpsi_1, \psi_2))\\
&=\CB((\VE, \VK_1, K_2),(\Vpsi, \Vpsi_1, \psi_2))\\
&\quad-\CB((\VE_H, \VK_{h,1}(\VE_H), K_{h,2}(\VE_H)), (\Vpsi, \Vpsi_1, \psi_2))\\
&=(\Vf, \Vpsi)-\CB((\VE_H, \VK_{h,1}(\VE_H), K_{h,2}(\VE_H)), (\Vpsi, \Vpsi_1, \psi_2))\\
&\quad-\CB_h((\VE_H, \VK_{h,1}(\VE_H), K_{h,2}(\VE_H)), (\Vpsi-\Vpsi_H, \Vpsi_1-\Vpsi_h, \psi_2-\psi_h))\\
&\quad+\CB_h((\VE_H, \VK_{h,1}(\VE_H), K_{h,2}(\VE_H)), (\Vpsi, \Vpsi_1, \psi_2))\\
&\quad-\underbrace{\CB_h(((\VE_H, \VK_{h,1}(\VE_H), K_{h,2}(\VE_H)), (\Vpsi_H, \Vpsi_h, \psi_h))}_{=(\Vf, \Vpsi_H)_{L^2}}.
\end{split}
\end{equation*}
With the definition of the residual \eqref{eq:residual}, this gives the following error identity for all $(\Vpsi, \Vpsi_1, \psi_2)\in \Hbf_0(\curl)\times L^2(\Om; \Hbf^1_{\sharp,0}(Y))\times L^2(\Om; H^1_{\sharp,0}(Y))$ and all $(\Vpsi_H, \Vpsi_h, \psi_h)\in \VV_{H, 0}^I\times L^2(\Om; \widetilde{\VW}_h^1(Y))\times L^2(\Om; \widetilde{W}_h^1(Y))$:
\begin{equation}
\label{eq:errorid}
\begin{split}
&\!\!\!\!\CB((e_0, e_1, e_2), (\Vpsi, \Vpsi_1, \psi_2))\\
&=-\langle\Res_h(\VE_H, \VK_{h,1}(\VE_H), K_{h,2}(\VE_H)), (\Vpsi-\Vpsi_H, \Vpsi_1-\Vpsi_h, \psi_2-\psi_h)\rangle\\
&\quad+(\CB_h-\CB)((\VE_H, \VK_{h,1}(\VE_H), K_{h,2}(\VE_H)), (\Vpsi, \Vpsi_1, \psi_2)).
\end{split}
\end{equation}
We choose $\Vpsi=e_0$, $\Vpsi_1=e_1$, $\psi_2=e_2$ and $\Vpsi_H=I_He_0$, $\Vpsi_h=I_he_1$, $\psi_h=I_he_2$ in the error identity \eqref{eq:errorid} with the interpolation operators $I_H$ from Lemma \ref{lem:schoeberlintpol} and $I_h$ from Lemma \ref{lem:clementintpol}. Using the coercivity of $\CB$, we obtain
\begin{equation*}
\begin{split}
\|(e_0&, e_1, e_2)\|^2_e\\
\leq &C (\underbrace{|\langle \Res_h(\VE_H, \VK_{h,1}(\VE_H), K_{h,2}(\VE_H)), (e_0-I_He_0, e_1-I_he_1, e_2-I_h e_2)\rangle|}_{:=\text{\Romannum{1}}}\\
&\qquad+\underbrace{|(\CB_h-\CB)((\VE_H, \VK_{h,1}(\VE_H), K_{h,2}(\VE_H)), (e_0, e_1, e_2))|}_{:=\text{\Romannum{2}}}).
\end{split}
\end{equation*}

To estimate \Romannum{1}, we insert the decomposition of $e_0-I_He_0=\nabla \theta+\Vz$ according to Lemma \ref{lem:schoeberlintpol} and thus obtain
\begin{equation*}
\begin{split}
\text{\Romannum{1}}&=\Bigl|\int_\Om\int_Y \mu^{-1}_h(\curl\VE_H+\curl_y\VK_{h,1})\cdot (\curl\Vz^*+\curl_y(e_1-I_he_1)^*)\\
&\qquad\qquad+\Div_y\VK_{h,1}\Div_y(e_1-I_he_1)^*\\
&\qquad\qquad-\kappa_h(\VE_H+\nabla_yK_{h,2})\cdot(\nabla\theta^*+\Vz^*+\nabla_y(e_2-I_he_2)^*)\, dydx\\
&\quad\;-\int_\Om\Vf(\nabla\theta^*+\Vz^*)\, dx\Bigr|.
\end{split}
\end{equation*}
Integrating by parts locally, inserting $\Vf_H-\Vf_H$ and noting $\Div\Vf=0$ yields
\begin{align*}
\text{\Romannum{1}}&= \Bigl|\sum_{j\in J}\int_{T_j}\sum_{i\in I}\int_{S_i}\curl_x(\mu_h^{-1}(\curl\VE_H+\curl_y\VK_{h,1}))\cdot \Vz^*\, dydx\\*
&\quad+\sum_{j \in J}\int_{\partial T_j}\sum_{i\in I}\int_{S_i}(\mu_h^{-1}(\curl\VE_H+\curl_y\VK_{h,1}))\times \Vn\cdot \Vz^*\, dyd\sigma\\
&\quad+\sum_{j\in J}\int_{T_j}\sum_{i\in I}\int_{S_i}\curl_y(\mu_h^{-1}(\curl\VE_H+\curl_y\VK_{h,1}))\cdot (e_1-I_he_1)^*\, dy dx\\*
&\quad+\sum_{j \in J}\int_{T_j}\sum_{i\in I}\int_{\partial S_i}(\mu_h^{-1}(\curl\VE_H+\curl_y\VK_{h,1}))\times \Vn\cdot (e_1-I_he_1)^*\, d\si dx\\
&\quad-\sum_{j\in J}\int_{T_j}\sum_{i\in I}\int_{S_i}\nabla_y(\Div_y \VK_{h,1})\cdot(e_1-I_he_1)^*\, dy dx\\*
&\quad+\sum_{j\in J}\int_{T_j}\sum_{i\in I}\int_{\partial S_i}\Div_y\VK_{h,1}\Vn\cdot(e_1-I_he_1)^*\, d\si dx\\
&\quad+\sum_{j\in J}\int_{T_j}\sum_{i\in I}\int_{S_i}\Div_x(\kappa_h(\VE_H+\nabla_yK_{h,2}))\,\theta^*\, dydx\\*
&\quad-\sum_{j\in J}\int_{\partial T_j}\sum_{i\in I}\int_{S_i}(\kappa_h(\VE_H+\nabla_yK_{h,2}))\cdot \Vn\;\theta^*\, dyd\si\\
&\quad+\sum_{j \in J}\int_{T_j}\sum_{i \in I}\int_{S_i}\Div_y(\kappa_h(\VE_H+\nabla_yK_{h,2}))\cdot(e_2-I_he_2)^*\, dy dx\\*
&\quad-\sum_{j \in J}\int_{T_j}\sum_{i \in I}\int_{\partial S_i}(\kappa_h(\VE_H+\nabla_yK_{h,2}))\cdot\Vn\, (e_2-I_he_2)^*\, d\si dx\\
&\quad-\sum_{j\in J}\int_{T_j}\sum_{i\in I}\int_{S_i}(\Vf_H+\Vf-\Vf_H+\kappa_h(\VE_H+\nabla_y K_{h,2}))\cdot \Vz^*\, dy dx\Bigr|.
\end{align*}
As $\mu^{-1}_h, \kappa_h$ are constant on the cells $T_j\times S_i$, the correctors $K_h$ are constant with respect to $x$ and linear with respect to $y$, and $\VE_H$ is linear with respect to $x$, all terms with two derivatives with respect to the same variable cancel out. We derive by a rearrangement of sums and the H\"older inequality
\begin{align*}
\!\!\!\!\text{\Romannum{1}}&= \Bigl|\sum_{(j,l)\in\CE(\CT_H)}\int_{F_{jl}}\Bigl[\int_Y(\mu^{-1}_h(\curl\VE_H+\curl_y\VK_{h,1}))\times \Vn\, dy\Bigr]_{F_{jl}}\cdot\Vz^*\, d\si\\
&\quad-\sum_{(j,l)\in \CE(\CT_H)}\int_{F_{jl}}\Bigl[\int_Y(\kappa_h(\VE_H+\nabla_yK_{h,2}))\cdot \Vn\, dy\Bigr]_{F_{jl}}\theta^*\, d\si\\
&\quad+\sum_{j\in J}\!\sum_{(i,k)\in \CE(\CT_h)}\!\int_{T_j}\int_{\tilde{F}_{ik}}[(\mu^{-1}_h(\curl\VE_H+\curl_y\VK_{h,1}))\times \Vn+\Div_y\VK_{h,1}\Vn]_{\tilde{F}_{ik}}\\*
&\qquad\qquad\qquad\qquad\qquad\quad\cdot(e_1-I_he_1)^*\, d\si dx\\
&\quad-\sum_{j\in J}\sum_{(i,k)\in \CE(\CT_h)}\int_{T_j}\int_{\tilde{F}_{ik}}\left[(\kappa_h(\VE_H+\nabla_yK_{h,2}))\cdot\Vn\right]_{\tilde{F}_{ik}}(e_2-I_he_2)^*\, d\si dx\\ 
&\quad+\sum_{j\in J}\int_{T_j}\Div_x\Bigl(\int_Y\kappa_h(\VE_H+\nabla_y K_{h,2})\, dy\Bigr)\theta^*\, dx\\
&\quad+\sum_{j\in J}\int_{T_j}\Bigl(\Vf-\Vf_H+\Vf_H+\int_Y \kappa_h(\VE_H+\nabla_y K_{h,2})\,dy\Bigr)\Vz^*\, dx\Bigr|\\
&\leq\sum_{(j,l)\in\CE(\CT_H)}\Bigl\|\Bigl[\int_Y(\mu^{-1}_h(\curl\VE_H+\curl_y\VK_{h,1}))\times \Vn\, dy\Bigr]_{F_{jl}}\Bigr\|\|\Vz\|_{L^2(F_{jl})}\\
&\quad+\sum_{(j,l)\in \CE(\CT_H)}\Bigl\|\Bigl[\int_Y(\kappa_h(\VE_H+\nabla_yK_{h,2}))\cdot \Vn\, dy\Bigr]_{F_{jl}}\Bigr\|\,\|\theta\|_{L^2(F_{jl})}\\
&\quad+\sum_{j\in J}\sum_{(i,k)\in \CE(\CT_h)}\Bigl\|\Bigl[(\mu^{-1}_h(\curl\VE_H+\curl_y\VK_{h,1}))\times \Vn+\Div_y\VK_{h,1}\Vn\Bigr]_{\tilde{F}_{ik}}\Bigr\|\\*
&\qquad\qquad\qquad\qquad\|(e_1-I_he_1)\|_{L^2(T_j\times\tilde{F}_{ik})}\\
&\quad+\sum_{j\in J}\sum_{(i,k)\in \CE(\CT_h)}\|[(\kappa_h(\VE_H+\nabla_yK_{h,2}))\cdot\Vn]_{\tilde{F}_{ik}}\|\,\|(e_2-I_he_2)\|_{L^2(T_j\times\tilde{F}_{ik})}\\ 
&\quad+\sum_{j\in J}\Bigl\|\Div_x\Bigl(\int_Y\kappa_h(\VE_H+\nabla_y K_{h,2})\, dy\Bigr)\Bigr\|\,\|\theta\|_{L^2(T_j)}\\
&\quad+\sum_{j\in J}\Bigl(\|\Vf-\Vf_H\|+\Bigl\|\Vf_H+\int_Y \kappa_h(\VE_H+\nabla_y K_{h,2})\,dy\Bigr\|\Bigr)\,\|\Vz\|_{L^2(T_j)}.
\end{align*}
Using the estimates for the interpolation operators from Lemmas \ref{lem:clementintpol} and \ref{lem:schoeberlintpol}, we get
\begin{align*}
\text{\Romannum{1}}&\leq \sum_{(j,l)\in\CE(\CT_H)}CH_{jl}^{1/2}\Bigl\|\Bigl[\int_Y(\mu^{-1}_h(\curl\VE_H+\curl_y\VK_{h,1}))\times \Vn\, dy\Bigr]_{F_{jl}}\Bigr\|\\*
&\qquad \qquad \qquad\qquad\;\,\|\curl e_0\|_{L^2(\tilde{\om}_{T_j})}\\
&\quad+\sum_{(j,l)\in \CE(\CT_H)}CH_{jl}^{1/2}\Bigl\|\Bigl[\int_Y(\kappa_h(\VE_H+\nabla_yK_{h,2}))\cdot \Vn\, dy\Bigr]_{F_{jl}}\Bigr\|\,\|e_0\|_{L^2(\tilde{\om}_{T_j})}\\
&\quad+\sum_{j\in J}\!\sum_{(i,k)\in \CE(\CT_h)}\!\!\!\!\!Ch_{ik}^{1/2}\Bigl\|[(\mu^{-1}_h(\curl\VE_H+\curl_y\VK_{h,1}))\!\times \!\Vn+\Div_y\VK_{h,1}\Vn]_{\tilde{F}_{ik}}\Bigr\|\\*
&\qquad\qquad\qquad\qquad\quad\;\,\|\nabla_y e_1\|_{L^2(T_j\times\om_{\tilde{F}_{ik}})}\\
&\quad+\sum_{j\in J}\sum_{(i,k)\in \CE(\CT_h)}Ch_{ik}^{1/2}\|[(\kappa_h(\VE_H+\nabla_yK_{h,2}))\cdot\Vn]_{\tilde{F}_{ik}}\|\,\|\nabla_ye_2\|_{L^2(T_j\times\om_{\tilde{F}_{ik}})}\\ 
&\quad+\sum_{j\in J}CH_j\Bigl\|\Div_x\Bigl(\int_Y\kappa_h(\VE_H+\nabla_y K_{h,2})\, dy\Bigr)\Bigr\|\|e_0\|_{L^2(\tilde{\om}_{T_j})}\\
&\quad+\sum_{j\in J}CH_j\Bigl(\|\Vf-\Vf_H\|+\Bigl\|\Vf_H+\int_Y \kappa_h(\VE_H+\nabla_y K_{h,2})\,dy\Bigr\|\Bigr)\\*
&\qquad\qquad\qquad\quad\|\curl e_0\|_{L^2(\tilde{\om}_{T_j})}.
\end{align*}
Applying the Cauchy-Schwarz inequality, we obtain the desired local estimators and terms like $\sum_{j\in J}\|\curl e_0\|_{L^2(\tilde{\om}_{T_j})}$. As the triangulation is shape regular, each element $T_j$ only appears in a finite number of these neighborhoods and this number can be bounded above by a uniform constant (independent of $H, h$). Thus, we derive
\begin{align*}
\text{\Romannum{1}}&\leq C\Bigl(\sum_{(j,l)\in\CE(\CT_H)}\eta_{jl,1}^2\Bigr)^{1/2}\,\|\curl e_0\|_{L^2(\Om)}+C\Bigl(\sum_{(j,l)\in \CE(\CT_H)}\eta_{jl,2}^2\Bigr)^{1/2}\,\|e_0\|_{L^2(\Om)}\\
&\quad+C\Bigl(\sum_{j\in J}\sum_{(i,k)\in \CE(\CT_h)}\eta_{j,ik,1}^2\Bigr)^{1/2}\,\|\nabla_y e_1\|_{L^2(\Om\times Y)}\\
&\quad+C\Bigl(\sum_{j\in J}\sum_{(i,k)\in \CE(\CT_h)}\eta_{j,ik,2}^2\Bigr)^{1/2}\,\|\nabla_y e_2\|_{L^2(\Om\times Y)}\\
&\quad+C\Bigl(\sum_{j\in J}\eta_{j,2}^2\Bigr)^{1/2}\,\|e_0\|_{L^2(\Om)}+C\Bigl(\sum_{j\in J}\zeta_j^2+\eta_{j,1}^2\Bigr)^{1/2}\,\|\curl e_0\|_{L^2(\Om)}.
\end{align*}
Of course, all norms of the errors $e_0$, $e_1$, and $e_2$ can simply be estimated by $\|(e_0, e_1, e_2)\|_e$.

To estimate \Romannum{2}, we just split the integral into local terms and use the H\"older inequality
\begin{align*}
\text{\Romannum{2}}&=\Bigl|\sum_{j\in J}\int_{T_j}\sum_{i\in I}\int_{S_i}(\mu^{-1}(x_j, y_i)-\mu^{-1}(x,y))(\curl \VE_H+\curl_y\VK_{h,1})\\*
&\qquad\qquad\qquad\qquad\cdot(\curl e_0^*+\curl_y e_1^*)\\*
&\qquad\qquad\qquad\;\,+(\kappa(x_j, y_i)-\kappa(x,y))(\VE_H+\nabla_y K_{h,2})\cdot(e_0^*+\nabla_y e_2^*)\, dydx\Bigr|\\*
&\leq\Bigl(\sum_{j\in J}\sum_{i \in I}\zeta_{ji}^2\Bigr)^{1/2}\,\|(e_0, e_1, e_2)\|_e.
\end{align*}
Dividing each estimate by $\|(e_0, e_1, e_2)\|_e$ and combining both gives us the a posteriori error estimate. 
\end{Proof}

For the proof of the lower bound we need local bubble functions. Let us denote by $\lambda_{T,l}$, $l=1,\ldots, 4$, the barycentric coordinates of a tetrahedron $T$ and by $\lambda_{F, l}$, $l=1,2,3$, the barycentric coordinates of a face $F$. The local bubble functions on elements and faces are defined as
$$\psi_T:=256\prod_{l=1}^4\lambda_{T,l}, \qquad \psi_F:= 27\prod_{l=1}^3\lambda_{F,l}.$$
They fulfill $0\leq\psi_T, \psi_F\leq 1$, $\supp\psi_T\subset T$, and $\supp\psi_F\subset\hat{\om}_F$.
We also define a continuation operator $P_F: L^\infty(F)\to L^\infty(\hat{\om}_F)$ as the constant extension of a function in the direction perpendicular to the face $F$, for details see \cite{Verf}. 
The following inequalities can be proven with standard scaling arguments and the properties of the bubble functions (see \cite[Proposition 3.37]{Verf} for details and the proof).

\begin{lemma}[Inverse inequalities]
\label{lem:bubbles}
For all $g\in \pz^k$ and all tetrahedra $T$ it holds
\begin{equation*}
\begin{split}
\|g\|_{L^2(T)}^2&\leq C |(g, \psi_T g)_T|,\\
\|\psi_T g\|_{L^2(T)}&\leq C\|g\|_{L^2(T)},\\
\|\nabla(\psi_T g)\|_{L^2(T)}&\leq C \diam(T)^{-1}\|g\|_{L^2(T)}.
\end{split}
\end{equation*}
Furthermore, for all $f\in \pz^k|_F$ and faces $F$ it holds
\begin{equation*}
\begin{split}
\|f\|_{L^2(F)}^2&\leq C |(g, \psi_F P_F(f))_F|,\\
\|\psi_F P_F(f)\|_{L^2(\om_F)}&\leq C\diam(F)^{1/2}\|f\|_{L^2(F)},\\
\|\nabla(\psi_F P_F(f))\|_{L^2(\om_F)}&\leq C \diam(F)^{-1/2}\|f\|_{L^2(F)}.
\end{split}
\end{equation*}
\end{lemma}

\begin{Proof}[Proof of Theorem \ref{thm:lowerbound}]
First of all, corresponding to the error terms we introduce the following functions:
\begin{align*}
\Vw_{j,1}(x)&=\psi_{T_j}(x)\Bigl(\Vf_H+\int_Y\kappa_h(\cdot, y)(\VE_H+\nabla_y K_{h,2}(\cdot, y))\,dy\Bigr)(x),\\
\Vw_{jl,1}(x)&=\psi_{F_{jl}}(x)P_{F_{jl}}\Bigl(\Bigl[\int_Y(\mu^{-1}_h(\curl\VE_H+\curl_y \VK_{h,1}))\times \Vn\, dy\Bigr]_{F_{jl}}\Bigr)(x),\\
\Vw_{j,ik,1}(x,y)&=\chi_{T_j}(x)\psi_{\tilde{F}_{ik}}(y)\\
&\quad P_{\tilde{F}_{ik}}([(\mu^{-1}_h(\curl\VE_H+\curl_y \VK_{h,1}))\times\Vn+\Div_y \VK_{h,1}\Vn]_{\tilde{F}_{ik}})(y),\\
w_{j,2}(x)&=\psi_{T_j}(x)\Div_x\Bigl(\int_Y \kappa_h(\cdot, y)(\VE_H+\nabla_yK_{h,2}(\cdot, y))\, dy\Bigr)(x),\\
w_{jl,2}(x)&=\psi_{F_{jl}}(x)P_{F_{jl}}\Bigl(\Bigl[\int_Y \kappa_h(\cdot, y)(\VE_H+\nabla_y K_{h,2}(\cdot, y))\cdot \Vn\, dy\Bigr]_{F_{jl}}\Bigr)(x),\\
w_{j,ik,2}(x,y)&=\chi_{T_j}(x)\psi_{\tilde{F}_{ik}}(y)P_{\tilde{F}_{ik}}([\kappa_h(\VE_H+\nabla_y K_{h,2})\cdot \Vn]_{\tilde{F}_{ik}})(y),
\end{align*}
where $\chi_A$ denotes the characteristic function of the set $A$.

Now the error indicators can be estimated separately: By partial integration they can be interpreted as the residual tested with the localized functions. The error identity and Lemma \ref{lem:bubbles} are then used each time to bound the indicators by the total error.

As $\VE_H$ is linear and $\VK_{h,1}$ constant with respect to $x$, we get with integration by parts and $\supp \Vw_{j,1}\subset T_j$
\begin{align*}
0&=\int_\Om\int_Y\curl_x(\mu^{-1}_h (\curl\VE_H+\curl_y \VK_{h,1}))\cdot \Vw_{j,1}^*\, dydx\\*
&=\int_\Om\int_Y\mu^{-1}_h(\curl\VE_H+\curl_y \VK_{h,1})\cdot \curl \Vw_{j,1}^*\, dydx.
\end{align*}
Therefore, we obtain for the first error indicator
\begin{align*}
\frac{\eta_{j,1}^2}{H_j^2}&\leq C\Bigl|\int_{T_j}\Bigl(\Vf_H+\int_Y\kappa_h(\VE_H+\nabla_y K_{h,2})\, dy\Bigr)\cdot \Vw_{j,1}^*\, dx\Bigr|\\*
&=C\left|-\langle\Res_h(\VE_H, \VK_{h,1}, K_{h,2}), (\Vw_{j,1}, 0, 0)\rangle+(\Vf_H-\Vf, \Vw_{j,1})_{L^2}\right|.
\end{align*}
If we choose $\Vpsi=\Vw_{j,1}$, $\Vpsi_1=0$, $\psi_2=0$ in the error identity \eqref{eq:errorid}, we get with Lemma \ref{lem:bubbles}
\begin{align*}
&\!\!\!\!|-\langle\Res_h(\VE_H, \VK_{h,1}, K_{h,2}), (\Vw_{j,1}, 0, 0)\rangle+(\Vf_H-\Vf, \Vw_{j,1})_{L^2}|\\*
&=|\CB((e_0, e_1, e_2), (\Vw_{j,1}, 0, 0))+(\CB-\CB_h)((\VE_H, \VK_{h,1}, K_{h,2}), (\Vw_{j,1}, 0, 0))\\*
&\quad+(\Vf_H-\Vf, \Vw_{j,1})_{L^2}|\\
&\leq C\|(e_0, e_1, e_2)\|_{e(T_j\times Y)}\|\Vw_{j,1}\|_{\Hbf(\curl, T_j)}+C\|\Vf_H-\Vf\|_{L^2(T_j)}\|\Vw_{j,1}\|_{L^2(T_j)}\\*
&\quad+C\Bigl(\sum_{i \in I}\zeta_{ji}^2\Bigr)^{1/2}\,\|\Vw_{j,1}\|_{\Hbf(\curl, T_j)}\\
&\leq C\|\Vf_H-\Vf\|_{L^2(T_j)}H_j^{-1}\, \eta_{j,1}+CH_j^{-2}\|(e_0,e_1,e_2)\|_{e(T_j\times Y)}\, \eta_{j,1}\\*
&\quad+CH_j^{-2}\Bigl(\sum_{i \in I}\zeta_{ji}^2\Bigr)^{1/2}\, \eta_{j,1}.
\end{align*}
All in all, after multiplying by $H_j^2\, \eta_{j,1}^{-1}$, this gives the local estimate for $\eta_{j,1}$.

For $\eta_{j,2}$ we get with the properties of $w_{j,2}$ and an integration by parts
\begin{align*}
\frac{\eta_{j,2}^2}{H_j^2}&\leq C \Bigl|\int_{T_j}\Div_x \Bigl(\int_Y \kappa_h(\VE_H+\nabla_y K_{h,2})\,dy\Bigr)w_{j,2}^*\, dx\Bigr|\\
&=C\Bigl|-\int_{T_j}\Bigl(\int_Y \kappa_h(\VE_H+\nabla_y K_{h,2})\, dy\Bigr)\cdot\nabla w_{j,2}^*\, dx\Bigr|\\
&=C\left|\langle\Res_h(\VE_H, \VK_{h,1}, K_{h,2}), (\nabla w_{j,2}, 0, 0)\rangle\right|,
\end{align*}
where the last equality holds because of $\curl\nabla=0$ and $(\Vf, \nabla w_{j,2})_{L^2}=0$. 
Now we can use the error identity \eqref{eq:errorid} with $\Vpsi=\nabla w_{j,2}$, $\Vpsi_1=0$, $\psi_2=0$ and obtain with Lemma \ref{lem:bubbles}
\begin{align*}
&\!\!\!\!|\langle\Res_h(\VE_H, \VK_{h,1}, K_{h,2}), (\nabla w_{j,2}, 0, 0)\rangle|\\
&=|-\CB((e_0,e_1, e_2), (\nabla w_{j,2}, 0, 0))+(\CB_h-\CB)((\VE_H, \VK_{h,1}, K_{h,2}), (\nabla w_{j,2}, 0, 0))|\\
&\leq C\Bigl(\|(e_0, e_1, e_2)\|_{e(T_j\times Y)}\|\nabla w_{j,2}\|_{L^2(T_j)}+\Bigl(\sum_{i \in I}\zeta_{ji}^2\Bigr)^{1/2}\,\|\nabla w_{j,2}\|_{L^2(T_j)}\Bigr)\\
&\leq C\Bigl(H_j^{-2}\|(e_0, e_1, e_2)\|_{e(T_j\times Y)}\, \eta_{j,2}+H_j^{-2}\Bigl(\sum_{i \in I}\zeta_{ji}^2\Bigr)^{1/2}\,\eta_{j,2}\Bigr).
\end{align*}
Multiplication by $H_j^2\, \eta_{j,2}^{-1}$ gives the local estimate for $\eta_{j,2}$.

For $\eta_{jl,1}$ we have the estimate
$$\frac{\eta_{jl,1}^2}{H_{jl}}\leq C\Bigl|\int_{F_{jl}}\Bigl[\int_Y(\mu^{-1}_h(\curl\VE_H+\curl_y \VK_{h,1}))\times \Vn\, dy\Bigr]_{F_{jl}}\cdot \Vw_{jl,1}^*\,d\si\Bigr|.$$
An integration by parts and the linearity of $\VE_H$ yield
\begin{align*}
&\int_{F_{jl}}\Bigl[\int_Y(\mu^{-1}_h(\curl\VE_H+\curl_y \VK_{h,1}))\times \Vn\, dy\Bigr]_{F_{jl}}\cdot \Vw_{jl,1}^*\,d\si\\*
&=\int_{\hat{\om}_{F_{jl}}}\int_Y\mu^{-1}_h(\curl\VE_H+\curl_y \VK_{h,1}) \cdot \curl \Vw_{jl,1}^*\, dydx\\
&=\langle\Res_h(\VE_H, \VK_{h,1}, K_{h,2}), (\Vw_{jl,1}, 0, 0)\rangle+\int_{\hat{\om}_{F_{jl}}}(\Vf-\Vf_H)\cdot \Vw_{jl,1}^*\, dx\\*
&\quad+\int_{\hat{\om}_{F_{jl}}}\Bigl(\Vf_H+\int_Y \kappa_h(\VE_H+\nabla_y K_{h,2})\, dy\Bigr)\cdot \Vw_{jl,1}^*\, dx.
\end{align*}
The integrals over $\hat{\om}_{F_{jl}}$ can be split into two element integrals over $T_j$ and $T_l$. In the third term we recognize the known error indicator $\eta_{j,1}$ and in the second term the data approximation error indicator $\zeta_j$. Using first the error identity \eqref{eq:errorid} with $\Vpsi=\Vw_{jl,1}$, $\Vpsi_1=0$, and $\psi_2=0$ and Lemma \ref{lem:bubbles}, we obtain
\begin{align*}
\frac{\eta_{jl,1}^2}{H_{jl}}&\leq C\Bigl(\|(e_0, e_1, e_2)\|_{e(\hat{\om}_{F_{jl}}\times Y)}\|\Vw_{jl,1}\|_{\Hbf(\curl, \hat{\om}_{F_{jl}})}\\*
&\qquad\;+\Bigl(\sum_i \zeta_{ji}^2+\zeta_{li}^2\Bigr)^{1/2}\,\|\Vw_{jl,1}\|_{\Hbf(\curl, \hat{\om}_{F_{jl}})}\\
&\qquad\; +H_j^{-1}\eta_{j,1}\|\Vw_{jl,1}\|_{L^2(T_j)}+H_l^{-1}\eta_{l,1}\|\Vw_{jl,1}\|_{L^2(T_l)}\\*
&\qquad\;+\|\Vf-\Vf_H\|_{L^2(T_j)}\|\Vw_{jl,1}\|_{L^2(T_j)}+\|\Vf-\Vf_H\|_{L^2(T_l)}\|\Vw_{jl,1}\|_{L^2(T_l)}\Bigr)\\
&\leq C\Bigl(\|(e_0, e_1, e_2)\|_{e(\hat{\om}_{F_{jl}}\times Y)}H_{jl}^{-1}\,\eta_{jl,1}+H_{jl}^{-1}\Bigl(\sum_i \zeta_{ji}^2+\zeta_{li}^2\Bigr)^{1/2}\,\eta_{jl,1}\\*
&\qquad\;+(H_j^{-1}\,\eta_{j,1}+H_l^{-1}\,\eta_{l,1})\, \eta_{jl,1}+(H_j^{-1}\,\zeta_j+H_l^{-1}\,\zeta_l)\, \eta_{jl,1}\Bigr).
\end{align*}
Due to the regularity of the triangulation the quotients $H_{jl}/H_j$ and $H_{jl}/H_l$ can be bounded above and below. Thus, multiplication by $H_{jl}\eta_{jl,1}^{-1}$ together with the already derived estimate for $\eta_{j,1}$ yields the desired estimate for $\eta_{jl,1}$.

For $\eta_{jl,2}$ we obtain by the properties of $w_{jl,2}$, an integration by parts, and $(\Vf, \nabla w_{jl,2})_{L^2}=0$
\begin{align*}
\frac{\eta_{jl,2}^2}{H_{jl}}&\leq C\Bigl|\int_{F_{jl}}\Bigl[\int_Y (\kappa_h(\VE_H+\nabla_y K_{h,2}))\cdot \Vn\, dy \Bigr]_{F_{jl}}w_{jl,2}^*\, d\si\Bigr|\\
&=C\Bigl|\int_{\hat{\om}_{F_{jl}}}\Div_x \Bigl(\int_Y\kappa_h(\VE_H+\nabla_y K_{h,2})\, dy\Bigr)w_{jl,2}^*\, dx\\*
&\qquad+\int_{\hat{\om}_{F_{jl}}}\int_Y \kappa_h(\VE_H+\nabla_y K_{h,2})\cdot \nabla w_{jl,2}^*\, dydx\Bigr|\\
&=C\Bigl|\int_{\hat{\om}_{F_{jl}}}\Div_x \Bigl(\int_Y\kappa_h(\VE_H+\nabla_y K_{h,2})\, dy\Bigr)w_{jl,2}^*\, dx\\*
&\qquad-\langle\Res_h(\VE_H, \VK_{h,1}, K_{h,2}), (\nabla w_{jl,2}, 0, 0)\rangle\Bigr|.
\end{align*}
When split into integrals over $T_j$ and $T_l$, the first term can be identified with $\eta_{j,2}$ and $\eta_{l,2}$, respectively. For the second term we can insert again the error identity \eqref{eq:errorid} with $\Vpsi=\nabla w_{jl,2}$, $\Vpsi_1=0=\psi_2$. With Lemma \ref{lem:bubbles} we then get
\begin{align*}
\frac{\eta_{jl,2}^2}{H_{jl}}&\leq C\bigl(\|(e_0, e_1, e_2)\|_{e(\hat{\om}_{F_{jl}}\times Y)}\|\nabla w_{jl,2}\|_{L^2(\hat{\om}_{F_{jl}})}+H_j^{-1}\,\eta_{j,2}\|w_{jl,2}\|_{L^2(T_j)}\\*
&\qquad+H_l^{-1}\,\eta_{l,1}\,\|w_{jl,2}\|_{L^2(T_l)} +\Bigl(\sum_i \zeta_{ji}^2+\zeta_{li}^2\Bigr)^{1/2}\,\|\nabla w_{jl,2}\|_{L^2(\hat{\om}_{F_{jl}})})\\
&\leq C\bigl(\|(e_0, e_1, e_2)\|_{e(\hat{\om}_{F_{jl}}\times Y)}H_{jl}^{-1}\, \eta_{jl,2}+\bigl(H_j^{-1}\,\eta_{j,2}+H_l^{-1}\,\eta_{l,2}\bigr)\, \eta_{jl,2}\\*
&\qquad +H_{jl}^{-1}\Bigl(\sum_i \zeta_{ji}^2+\zeta_{li}^2\Bigr)^{1/2}\, \eta_{jl,2}\bigr).
\end{align*} 
Together with the already derived estimate for $\eta_{j,2}$ this gives us the local estimate for $\eta_{jl,2}$.

For $\eta_{j,ik,1}$ we have the estimate
\begin{align*}
\frac{\eta_{j,ik,1}^2}{h_{ik}}&\leq C\Bigl|\int_{T_j}\int_{\tilde{F}_{ik}}[(\mu^{-1}_h(\curl \VE_H+\curl_y \VK_{h,1}))\times \Vn+\Div_y \VK_{h,1}\Vn]_{\tilde{F}_{ik}}\\*
&\qquad\qquad\qquad\cdot \Vw_{j,ik,1}^*\, d\si dx\Bigr|.
\end{align*}
With an integration by parts and the linearity of $\VK_{h,1}$ with respect to $y$ we obtain
\begin{align*}
&\!\!\int_{T_j}\int_{\tilde{F}_{ik}}[(\mu^{-1}_h(\curl \VE_H+\curl_y \VK_{h,1}))\times \Vn+\Div_y \VK_{h,1}\Vn]_{\tilde{F}_{ik}}\cdot \Vw_{j,ik,1}^*\, d\si dx\\
&=\int_{T_j}\int_{\hat{\om}_{\tilde{F}_{ik}}}\!\!\mu^{-1}_h(\curl\VE_H+\curl_y\VK_{h,1})\!\cdot\! \curl_y \Vw_{j,ik,1}^*\!+\!\Div_y \VK_{h,1}\Div_y \Vw_{j,ik,1}^*\, dydx\\
&=\langle\Res_h(\VE_H, \VK_{h,1}, K_{h,2}), (0, \Vw_{j,ik,1}, 0)\rangle.
\end{align*}
Inserting the error identity \eqref{eq:errorid} with $\Vpsi=0=\psi_2$, $\Vpsi_1=\Vw_{j,ik,1}$ and using Lemma \ref{lem:bubbles}, we obtain
\begin{align*}
\frac{\eta_{j,ik,1}^2}{h_{ik}}&\leq C\bigl(\|(e_0, e_1, e_2)\|_{e(T_j\times \hat{\om}_{\tilde{F}_{ik}})}\|\nabla \Vw_{j,ik,1}\|_{L^2(T_j\times \om_{\tilde{F}_{ik}})}\\*
&\qquad+(\zeta_{ji}^2+\zeta_{jk}^2)^{1/2}\|\nabla \Vw_{j,ik,1}\|_{L^2(T_j\times \hat{\om}_{\tilde{F}_{ik}})})\\*
&\leq C(h_{ik}^{-1}\,(\zeta_{ji}^2+\zeta_{jk}^2)^{1/2}\, \eta_{j,ik,1}+h_{ik}^{-1}\,\|(e_0, e_1, e_2)\|_{e(T_j\times \hat{\om}_{\tilde{F}_{ik}})}\, \eta_{j,ik,1}\bigr).
\end{align*}

For $\eta_{j,ik,2}$ we derive with integration by parts and the linearity of $K_{h,2}$ with respect to $y$
\begin{align*}
\frac{\eta_{j,ik,2}^2}{h_{ik}}&\leq C\Bigl|\int_{T_j}\int_{\tilde{F}_{ik}}[(\kappa_h(\VE_H+\nabla_y K_{h,2}))\cdot \Vn]_{\tilde{F}_{ik}}\cdot w_{j,ik,2}^*\, d\si dx\Bigr|\\*
&=C\Bigl|\int_{T_j}\int_{\hat{\om}_{\tilde{F}_{ik}}}\kappa_h(\VE_H+\nabla_y K_{h,2})\cdot \nabla_yw_{j,ik,2}^*\, dydx\Bigr|\\*
&=C\bigl|-\langle\Res_h(\VE_H, \VK_{h,1}, K_{h,2}), (0,0,w_{j,ik,2})\rangle\bigr|.
\end{align*}
Inserting once more the error identity \eqref{eq:errorid}, this time with $\Vpsi=0$, $\Vpsi_1=0$, and $\psi_2=w_{j,ik,2}$, we obtain with Lemma \ref{lem:bubbles}
\begin{align*}
\frac{\eta_{j,ik,2}^2}{h_{ik}}&\leq C\bigl(\|(e_0,e_1, e_2)\|_{e(T_j\times \hat{\om}_{\tilde{F}_{ik}})}\|\nabla_y w_{j,ik,2}\|_{L^2(T_j\times \hat{\om}_{\tilde{F}_{ik}})}\\*
&\qquad+(\zeta_{ji}^2+\zeta_{jk}^2)^{1/2}\,\|\nabla_y w_{j,ik,2}\|_{L^2(T_j\times \hat{\om}_{\tilde{F}_{ik}})}\bigr)\\
&\leq C\bigl(h_{ik}^{-1}\,\|(e_0,e_1, e_2)\|_{e(T_j\times \hat{\om}_{\tilde{F}_{ik}})}\, \eta_{j,ik,2}+h_{ik}^{-1}\,(\zeta_{ji}^2+\zeta_{jk}^2)^{1/2}\, \eta_{j,ik,2}\bigr).
\end{align*}
This gives us the local estimate for $\eta_{j,ik,2}$.
The global estimate now follows by summing up the local estimates. 
\end{Proof}

\section*{Conclusion} 
In this paper, we suggested a new Heterogeneous Multiscale Method (HMM) for the time-harmonic Maxwell equations. The basis is a homogenization result for a curl-curl-problem obtained with two-scale convergence. Divergence-regularization is applied to the corrector for the curl and thus we can get rid of the divergence-free constraint. The HMM can bee seen as direct finite element discretization with numerical quadrature of the two-scale homogenized equation, which is the crucial observation for the numerical analysis.   
The a priori error analysis shows that the HMM converges with the same rates as finite element methods for curl-curl-problems without oscillating coefficients (Theorems \ref{thm:aprioriH1}, \ref{thm:aprioridual}). The a posteriori error estimators are reliable and efficient (Theorems \ref{thm:aposteriori}, \ref{thm:lowerbound}) and can be used for adaptive algorithms future work.

\end{document}